\documentclass[a4paper,11pt]{amsart}

\usepackage[T1]{fontenc}
\usepackage[utf8]{inputenc}
\usepackage[left=2cm,right=2cm,top=2cm,bottom=2cm]{geometry}

\usepackage{mathtools}
\usepackage{amssymb}
\usepackage{dsfont}

\usepackage{enumitem}
\usepackage{graphicx}
\usepackage{subcaption}

\usepackage[square,sort,comma,numbers]{natbib}
\usepackage[colorlinks=true,linkcolor=blue,citecolor=blue,urlcolor=blue,breaklinks]{hyperref}
\usepackage{cleveref}

\numberwithin{equation}{section}
\numberwithin{figure}{section}

\newcommand{\M}{\mathcal{M}}
\newcommand{\cL}{\mathcal{L}}
\newcommand{\cC}{\mathcal{C}}
\newcommand{\cT}{\mathcal{T}}

\newcommand{\sfL}{\mathsf{L}}
\newcommand{\sfT}{\mathsf{T}}

\newcommand{\dt}{\partial_t}
\newcommand{\R}{\mathbb{R}}
\newcommand{\Rd}{\R^d}
\newcommand{\RRd}{\R^d\times\R^d}
\newcommand{\gradv}{\nabla_v}
\newcommand{\gradx}{\nabla_x}

\renewcommand{\d}{\,\mathrm{d}}
\newcommand{\dv}{\,\mathrm{d}v}
\newcommand{\dx}{\,\mathrm{d}x}
\def\1{\mathds{1}}

\newcommand{\bangle}[1]{\langle #1\rangle}
\newcommand{\wangle}[1]{\lfloor #1 \rceil}
\DeclarePairedDelimiter{\abs}{\lvert}{\rvert}
\DeclarePairedDelimiter{\norm}{\lVert}{\rVert}

\let\eps\varepsilon

\theoremstyle{plain}
\newtheorem{lemma}{Lemma}[section]
\newtheorem{thm}[lemma]{Theorem}
\newtheorem{prop}[lemma]{Proposition}
\newtheorem{hypothesis}{Hypothesis}

\title[Kinetic Fokker-Planck equations with (very) weak velocity confinements]{Convergence to a non-explicit steady state in non-factorized kinetic Fokker-Planck equations with (very) weak velocity confinements}

\author{\'Emeric Bouin}
\address[E. Bouin]{CEREMADE - Université Paris-Dauphine, PSL Research University, UMR CNRS 7534, Place du Mar\'echal de Lattre de Tassigny, 75775 Paris Cedex 16, France.}
\email{bouin@ceremade.dauphine.fr}

\author{Luca Ziviani}
\address[L. Ziviani]{CEREMADE - Université Paris-Dauphine, PSL Research University, UMR CNRS 7534, Place du Mar\'echal de Lattre de Tassigny, 75775 Paris Cedex 16, France.}
\email{ziviani@ceremade.dauphine.fr}

\date{\today}

\begin{document}

    \begin{abstract}
    In this article, we prove several convergence results for kinetic Fokker-Planck equations with strong space confinement but local equilibria with sub-Gaussian or polynomial tails and non-explicit global steady states. We extend the results of \cite{C21} to a wider class of fat-tailed local equilibria, with rates of convergence in a large class of weighted $\sfL^1$ spaces. We complement our results with numerical simulations to investigate the shape of the non-explicit steady state. 
    \end{abstract}
 %   \subjclass[2010]{60J60,35Q84,82C40,35B27,60K50,60G52,76P05}
\maketitle
    
\tableofcontents

\section{Introduction}
This work is devoted to the kinetic Fokker-Planck equation,
\begin{equation}\label{L}
    \partial_t f=\mathcal{L}f:=-v\cdot\gradx f+\gradx V\cdot \gradv f+\gradv\cdot \left(\M\gradv\left(\frac{f}{\M}\right)\right) 
\end{equation}
where the unknown function $f=f(t,x,v)$ depends on time $t\in[0,\infty)$, position $x\in\Rd$ and velocity $v\in\Rd$. We consider external potentials $V\colon \Rd\to\R$ of the form
\begin{equation}\label{eq:V}
    V(x)=\frac{\wangle{x}^\alpha}{\alpha}, \qquad \forall x\in\Rd,
\end{equation}
with $\alpha>1$, where we denote $\wangle{\cdot}=\sqrt{1+\abs{\cdot}^2}$. The local equilibrium may take one of the following forms. Either
\begin{align}\label{eq:Mexp}
     \M(v)=c_\beta^{-1} e^{-\frac{\wangle{v}^\beta}{\beta}}, \qquad \forall v\in\Rd,
\end{align}
with $\beta>0$ and $c_\beta = \int_{\R^d} \exp\left(-\frac{\wangle{v}^\beta}{\beta}\right)\dv$, or
\begin{align}\label{eq:Mpoly}
     \M(v)=d_\gamma^{-1} \wangle{v}^{-d-\gamma}, \qquad \forall v\in\Rd, 
\end{align}
with $\gamma>0$ and $d_\gamma =\int_{\R^d} \wangle{v}^{-d-\gamma}\dv$. The evolution equation is supplemented with an initial datum $f_0\in\sfL^1$. The equation is mass conservative, that is,
\[
\int_{\R^d \times \R^d} f\dx\dv = \int_{\R^d \times \R^d} f_0\dx\dv,\qquad\text{for all } t\geq 0.
\]

Equation \eqref{L} is usually used to describe the evolution in time of the density function of particles under the influence of an external macroscopic potential $V$. It is a very suitable physical model for gases, plasmas, or stellar systems, but it has also numerous applications in biology, chemistry, and material  science. The operator $\cL$ in \eqref{L} can be split into $\cL=\sfL - \sfT$, where $\sfT$ is a transport operator and $\sfL$ is a collision operator. The transport operator $\sfT$ is 
\[
\sfT f = v\cdot\gradx f-\gradx V\cdot \gradv f,
\]
and corresponds to the Lie derivative of $f$ along the Hamiltonian vector field associated with the energy
\begin{equation}
    E(x,v)=\frac{\abs{v}^2}{2} + V(x).
\end{equation}
The collision operator $\sfL$ is the Fokker-Planck operator,
\[
\sfL f =\gradv\cdot \left( \M\gradv\left(\frac{f}{ \M}\right)\right) = \Delta_v f -\gradv\cdot\left(\frac{\gradv \M}{\M}f\right)
\]
Observe that our assumptions on $\M$ (\eqref{eq:Mexp} or \eqref{eq:Mpoly}) lead to power law drifts, that is 
\begin{equation*}
    \frac{\gradv \M}{\M} = -\wangle{v}^{\beta-2}v
\end{equation*}
if $\M$ is given by \eqref{eq:Mexp}, and
\begin{equation*}
    \frac{\gradv \M}{\M} =-(d+\gamma) \wangle{v}^{-2}v
\end{equation*}
if $\M$ is given by \eqref{eq:Mpoly}. Motivated by these expressions, it is natural to define the associated parameter $\beta = 0$ when $\M$ is given by \eqref{eq:Mpoly}. This choice will turn out to be particularly useful to simplify the notation along the proof.

In the case $\beta=2$, equation \eqref{L} turns out to be the classical kinetic Fokker-Planck equation, with an explicit global Gibbs state,
\begin{equation}\label{eq:Gibbs}
G = Z^{-1}\exp(-E)    
\end{equation}
for some normalisation constant $Z>0$. This equation has been extensively studied and various hypocoercivity methods have been used \cite{Dolbeault2015,MR2576899,MR2813582,Helffer2005,Hrau2004,MR2294477,Villani2009} to show exponential convergence to the steady state in weighted $\sfL^2$ norms. Still with $\beta=2$, in \cite{Cao2019}, the author proves convergence for weaker potential \eqref{eq:V} with $\alpha\in(0,1)$.
For other results with weak confinement we refer the interested reader to \cite{Bakry2008,MR2499863}.

On the other hand, when the local equilibrium $\M$ is not a Gaussian, the transport operator and the collision operator do not share the same kernel. This fact makes the existence of a stationary state $G$ less trivial: indeed such a function must solve the equation $\sfT G=\sfL G$ where both sides of the equation are not zero. It is therefore necessary to use non-explicit methods, like the Harris theorem, to show the existence of such an equilibrium.

For further use, we define the formal dual operator $\cL^*$ by
\begin{equation}\label{Lstar}
    \mathcal{L}^*m:=v\cdot\gradx m-\gradx V\cdot \gradv m +\Delta_v m+\frac{\gradv  \M}{ \M}\cdot\gradv m  
\end{equation}
For a positive weight function $m\colon \RRd\to \R$, we denote by $\sfL^1(m)$ the functional space defined by the norm
\[
\norm{f}_{\sfL^1(m)}:=\iint_{\R^d \times \R^d} \abs{f(x,v)}m(x,v)\dx\dv,
\]
and we denote by $\sfL^1$ the standard unweighted space.

The main results of this paper are existence and uniqueness of a steady state $G$ to \eqref{L} and convergence of solutions towards $G$. The first main result concerns the kinetic Fokker-Planck equation when $ \M$ has thin tails as in \eqref{eq:Mexp}.

\begin{thm}\label{thm:mainEXP}
Let $V$ be given by \eqref{eq:V} with $\alpha>1$ and $\M$ by \eqref{eq:Mexp} with $\beta > 0$. Then,
\begin{enumerate}
    \item there exists a unique positive normalised steady state $G\in \sfL^1\left(\exp({\delta E^{\frac{\min\{\beta,2\}}{2}}})\right)$ with $\delta>0$ small enough,
    \item for every $\theta\in(0, \min(1,\frac{\beta}{2})]$, there exists $\lambda>0$ such that, for any $f_0\in\sfL^1(e^{\delta E^\theta})$, the associated solution $f$ to \eqref{L} satisfies 
\begin{equation*}
    \left\Vert f(t)- \left(\int_{\R^d \times \R^d} f_0\right)G\right\Vert_{\sfL^1}\lesssim \footnote{We use the notation $\mathsf a\lesssim \mathsf b$ if there exists a constant $\mathsf c>0$ such that $\mathsf a\leq \mathsf{c}\,\mathsf{b}$.}
    e^{-\lambda t^\theta} \left\Vert f_0- \left(\int_{\R^d \times \R^d} f_0 \right) G\right\Vert_{\sfL^1(e^{\delta E^\theta})}.
\end{equation*}
\end{enumerate}

\end{thm} 
This theorem extends the results of \cite{C21}, where only $\beta\geq2$ is covered. We provide the analogous result for the missing range of parameters $0<\beta<2$, and observe a slower decay rate. Notice that, for $\beta\geq 2$, we recover the exponential decay, as expected. Our proof can actually cover the result of \cite{C21}.

For fat-tailed local equilibrium \eqref{eq:Mpoly}, we have the following result.
\begin{thm}\label{thm:mainPOLY}
Let $V$ be given by \eqref{eq:V} with $\alpha>1$ and $\M$ by \eqref{eq:Mpoly} with $\gamma > 0$. Then,
\begin{enumerate}
    \item there exists a unique positive normalised steady state
    \[
    G\in \bigcap_{0<s<\frac\gamma2} \sfL^1(E^s)
    \]
    \item for every $k\in(1, 1+\frac\gamma2)$ and for any $f_0\in\sfL^1(E^k)$, the associated solution $f$ to \eqref{L} satisfies 
    \begin{equation*}
    \left\Vert f(t)- \left(\int_{\R^d \times \R^d} f_0 \right) G\right\Vert_{\sfL^1}\lesssim \frac{1}{\bangle{t}^{k}}\left\Vert f_0-\left(\int_{\R^d \times \R^d} f_0 \right)G\right\Vert_{\sfL^1(E^k)}.
\end{equation*}
\end{enumerate}
\end{thm}
Actually, the latter convergence rate result is also correct when $\M$ has the form \eqref{eq:Mexp}, and for any $k>1$. We shall mention it as a side result. 

\begin{prop}
Let $V$ be given by \eqref{eq:V} with $\alpha>1$ and $\M$ by \eqref{eq:Mexp} with $\beta > 0$. Let $k>1$, let $f$ be a solution to \eqref{L}, with initial data $f_0 \in \sfL^1(E^k)$ and let $G$ be the unique steady state of \Cref{thm:mainEXP}. For every $k>1$ fixed, we have
    \begin{equation}\label{eq:HPE}
        \left\Vert f(t)- \left(\int_{\R^d \times \R^d} f_0 \right)G \right\Vert_{\sfL^1}\lesssim\frac{1}{\bangle{t}^{k}}\left\Vert f_0- \left(\int_{\R^d \times \R^d} f_0 \right)G\right\Vert_{\sfL^1(E^k)}.
    \end{equation}   
\end{prop}
    
It is worth noticing that
\Cref{thm:mainEXP} and \Cref{thm:mainPOLY} are true even if the potential $V$ satisfies the bounds
\[
V(x)\asymp\footnote{We write $\mathsf a\asymp \mathsf b$
if both $\mathsf a\lesssim \mathsf b$ and $\mathsf b\lesssim \mathsf a$ hold.} \wangle{x}^\alpha
\]
and
\[
x\cdot\gradx V\asymp \abs{x}^2\wangle{x}^{\alpha-2}
\]
for all $x\in B_R^c$, where $B_R$ denotes the ball centered at the origin with radius $R>0$. Regarding the variable $v$ one could also lighten the exact $\tfrac{\gradv\M}{\M}=-\wangle{v}^{\beta-2}v$ in the same flavour for all the results above, but at the expense of heavier computations. 

Let us comment qualitatively the results of the main theorems. As we said, the transport operator and the collision operator do not always share the same kernel, as $\ker \sfT$ is composed by states which can be written as a function of the energy $E$, while $\ker \sfL$ is spanned by $\M$. Only the particular case when $\M$ is as in \eqref{eq:Mexp} with $\beta = 2$, the Gibbs state \eqref{eq:Gibbs} is the unique and explicit steady state of the equation. For the remaining cases, the expression of the steady state is not explicit and it would be extremely interesting to get some insights about its shape, in particular the decay in $x$ and $v$. In \Cref{thm:mainEXP}, with $\beta\in(0,2)$, we found the steady state $G$ in the space $\sfL^1\left(\exp({\delta E^{\frac{\beta}{2}}})\right)$; this might suggest that, at least in a first order approximation, the steady state is decaying like a function of the energy. Moreover, due to the strength of the weight, the decay has to be sufficiently fast, for example it can be of the form
\[
G\approx e^{-\delta' E^{\frac{\beta}{2}}},%\sim e^{-\frac{\wangle{v}^\beta}{\beta}}
\]
for some $\delta'>0$. It is reasonable to believe that this might actually be a good approximation of $G$ if we observe the decay in the variable $v$ individually. Indeed, we can see that, as $\abs{v}\to \infty$ and for suitable $\delta'$, the above approximation decays in $v$ like $\M$.

The same idea may apply to \Cref{thm:mainPOLY}, this time the collision operator is pushing towards a polynomially decaying function as in \eqref{eq:Mpoly}. We expect that the steady state might resemble a function of the energy and that, being in $\bigcap_{0<s<\frac\gamma2} \sfL^1(E^s)$, such approximation has energy moments of all orders strictly smaller than $\frac{\gamma}{2}$. Hence, we might have
\[
G\approx E^{-\frac{d+\gamma}{2}}.
\]
This implies in particular that $G$ has all moments in $v$ strictly under the order $\gamma$, exactly as $\M$. 

This work can be considered as an extension of \cite{C21}, where the kinetic Fokker-Planck equation \eqref{L} with local equilibrium \eqref{eq:Mexp} with $\beta\geq 2$ is considered. In that work, Cao first proved the existence of a unique steady state and convergence in weighted $\sfL^1$ norms by using the Harris theorem. In a second step, this result is extended to $\sfL^p$ spaces thanks to regularising properties of the semigroup $S_\cL$ and extension arguments. In the particular cases where the steady state does not have an explicit expression, Lyapunov functions play an important role, because they allow to deduce the existence of stationary states and a suitable weighted $\sfL^1$ space. However, finding sharp Lyapunov functions turns out to be particularly hard, especially for weaker confinements in velocity. We aim to complete the picture left open by \cite{C21}, this is one of the main contributions of the present paper. Compared to that work, we improved the expression of the Lyapunov functions with sharper terms in order to extend the convergence result to fatter tailed $\M$. With more precise estimates, we are able to cover the remaining cases, when $\M$ is given in \eqref{eq:Mexp} for $0<\beta<2$ and \eqref{eq:Mpoly} for $\gamma>0$. However, unlike \cite{C21}, the proof does not use the geometric version of Harris's theorem, but rather its sub-geometric version, which relies on a weak Lyapunov condition.

We shall stress that the condition $\alpha > 1$ is unavoidable in our computations: our main concern is the strength of the velocity drift. Few results are actually known when $\alpha \leq 1$: only the very special case $\beta = 2$ is discussed in \cite{Cao2019} and the arguments are there based on the fact that the steady state $G$ is known and explicit. The origin of the difficulty in other situations is that weak Lyapunov functions are not known at the moment, and results in this direction are the aim of a subsequent work. Some results when the collision operator is of scattering type are given in \cite{CCEY20}, but do not cover Fokker-Planck operators.

The main tool we will use is the Harris theorem, which is nowadays very popular to prove the existence of a stationary measure for stochastic semigroups and convergence toward it. The first result on the stability of Markov processes goes back to Doeblin \cite{Doeblin1940}, who showed exponential convergence to a stationary measure assuming that the transition probabilities possess a uniform lower bound. In \cite{H56}, some sufficient conditions for the existence of a stationary measure are discussed and, in the last few years, these results have been further developed. In particular, we mention \cite{HM11,DMT95,MT09,MT93,CM21,H16} for such theorems, which are now known as Harris-type results. In this work, we need the sub-geometric version of it, which we will recall in \Cref{sec:Harris} for self-consistency. For more information about the sub-geometric Harris theorem, see \cite{MR2499863,H16,CM21}. In the last decades, the Harris theorem has been successfully applied to many kinetic equations: such as the kinetic Fokker-Planck equation \cite{C21}, the homogeneous Fokker-Planck equation \cite{Kavian2021}, the linear Boltzmann equation \cite{CCEY20,EM24} and the run-and-tumble equation in biology \cite{BEZ25,EY24,EY23}. See also the review \cite{Yolda2023}. Two main ingredients are needed for the Harris theorem: a Lyapunov condition and a positivity condition. In the case of Fokker-Planck equations, the Lyapunov condition can be verified directly using the adjoint operator $\cL^*$, on the other hand, the positivity property corresponds to proving a certain Harnack inequality.

Harnack inequalities for kinetic Fokker-Planck equation were first proved in \cite{MR3923847} by a
non-constructive argument, later was revised by \cite{GM22} with a new constructive De Giorgi approach based on ideas of \cite{PP04}. See also \cite{Guerand_Imbert_2023} for another constructive proof by the Moser-Kružkov approach. As the name suggests, these methods are adaptations from the De Giorgi argument \cite{DG56,DG57} for proving H\"older regularity for elliptic equations. Moreover, Harnack inequalities have also been obtained for the long-range Boltzmann equation in \cite{MR4049224} and with applications to the Landau equation \cite{MR3923847}.

It is worth mentioning that a very similar model to \eqref{L} has been studied in \cite{BouinEmericandDolbeault2024}, replacing the transport operator $\sfT$ with another Hamiltonian operator more compatible with the collision operator. This greatly simplifies the problem as the steady state is explicit, so no Harris type arguments are needed. The authors managed to perform an $\sfL^2$-hypocoercivity method type method, taking its foundations in \cite{Dolbeault2015}. As in the present paper, slower convergence rates were obtained due to the heavy tails of $\M$.  

Mathematical models featuring non-Gaussian local equilibria have several application in modeling complex physical and socio-economic systems. In quantitative finance and wealth distribution dynamics, for instance, empirical distributions systematically exhibit fat-tailed behavior \cite{BCG26}. Within this framework, kinetic Fokker–Planck equations frequently serve as continuous approximations of the underlying Boltzmann equations, particularly in microscopic models governing wealth exchange and opinion dynamics \cite{Cordier2005,T06}.

The rest of the paper is organised as follows. In \Cref{sec:Lyapunov}, we provide many weak Lyapunov functions, first for the case \eqref{eq:Mexp} and then for \eqref{eq:Mpoly}. In \Cref{sec:positivity}, we prove the positivity condition required by the Harris theorem. Unlike in \cite{C21}, we do not show a regularisation result of the type $\sfL^1\to \cC^\infty$ about the semigroup $S_\cL$, but we use most recent results on Harnack inequalities as in \cite{CGMM24}. Finally, in \Cref{sec:proofs} we prove the main results. We end the paper by recalling the Harris theorem in \Cref{sec:Harris}.

\section{Weak Lyapunov functions}\label{sec:Lyapunov}
In this section, we seek weight functions $m \colon \R^d \times \R^d \to[1, \infty)$ verifying the so-called Lyapunov condition, that is
\[
\mathcal{L}^*m\leq C\1_{B_R} -\phi(m),
\]
where $C,\,R > 0$ and $\phi\colon\R_+\to \R_+$ is a concave function such that $\lim_{\abs{(x,v)}\to\infty}m(x,v)=\infty$. We provide below a family of such weights, in both cases $\M$ defined as in \eqref{eq:Mexp} and in \eqref{eq:Mpoly}. These functions are crucial to apply the Harris theorem of \Cref{sec:Harris}. 

We begin by introducing a useful function $H\colon\RRd\to \R$, defined by
\begin{equation}\label{eq:H}
    H(x,v) = E^\ell+\eps \, \frac{\wangle{x}^{A}}{\wangle{v}^{B}} \, x\cdot v,
\end{equation}
for some $\ell \geq 2$ , $\eps \in(0,1)$ and real constants $A \geq 0$ and $B  \in (0,1)$ to be adapted later on. 

\subsection{Preliminaries}
The function $H$ will be used to construct the Lyapunov functions $m$, so, to make the following calculations clearer, in the next lemma we provide some useful estimates about it. We highlight that the explicit constants are important for the rest of the paper.
\begin{lemma}\label{lem:tech}
Let $\ell\geq 2$, $\eps\in(0,1)$ and consider the function $H$ defined in \eqref{eq:H} with $A,B$ satisfying 
\begin{align}\label{eq:conAB}
    0< A &< \alpha(\ell-1) -1, & 0< B< 1.
\end{align}
Then there exists a constant $c>0$ such that, for $\eps>0$ small enough, $H$ is positive and satisfies the following bounds
\begin{align}\label{eq:boundH}
    (1- c \eps)E^\ell \leq H\leq (1+ c \eps)E^\ell
\end{align}
and
\begin{align}\label{eq:gradH}
    (1-\eps) \ell^2 E^{2\ell-2}\abs{v}^2\leq &\left\vert \nabla_v H \right\vert^2 \lesssim E^{2\ell-2}\wangle{v}^{2}. 
\end{align}
Moreover, there exists a constant $C>0$ such that 
\begin{equation}\label{eq:L*H}
    \begin{split}
        \cL^*H  \leq & \left[ (\ell-1) \frac{\abs{v}^2}{E} + d + v\cdot\frac{\gradv \M}{\M}\right]\ell E^{\ell-1} - (1-B)\eps\frac{\wangle{x}^{A + \alpha -2}\abs{x}^2}{\wangle{v}^{B}}\\
        &\qquad+ C\eps\frac{\wangle{x}^{A}}{\wangle{v}^{B}}\left[ \abs{v}^2 + \wangle{v}^{(\beta -1)_+}  \abs{x}\right]
    \end{split}
\end{equation}
and for any $\Phi\in\cC^2(\R_+)$, 
\begin{equation}\label{eq:L*Phi(H)}
   \mathcal{L}^*(\Phi(H)) = \Phi'(H) \mathcal{L}^*(H) + \Phi''(H)\abs{\gradv H}^2.
\end{equation}
\end{lemma}
\begin{proof}[{\bf Proof of \Cref{lem:tech}}]

\noindent\textbf{\# Proof of \eqref{eq:boundH}}. Thanks to the Cauchy-Schwarz inequality and \eqref{eq:conAB}, we bound the mixed term as 
\begin{align*}
    \frac{\wangle{x}^{A}}{\wangle{v}^{B}}\abs*{x\cdot v}\leq \wangle{x}^{1+A}\wangle{v}^{1-B}  \leq \wangle{x}^{\alpha(\ell-1)} \leq  \alpha^{\frac{\ell-1}{\alpha}} E^{\ell-1}  \leq c E^{\ell}
\end{align*}
with $c = 2\alpha^{\frac{\ell-1}{\alpha}}$. Therefore, by definition of $H$, we conclude that \eqref{eq:boundH} is satisfied.

\noindent\textbf{\# Proof of \eqref{eq:gradH}}. By a direct computation, we have that
\begin{align*}
\nabla_v H  = J_1 + \eps J_2.
\end{align*}
where
\begin{align*}
    J_1& = \ell E^{\ell - 1} v & J_2&=  \frac{\wangle{x}^{A}}{\wangle{v}^{B}} \, x - B \frac{\wangle{x}^{A} (x\cdot v) }{\wangle{v}^{B+2}} \, v 
\end{align*}
With these notations, we have $\abs{\gradv H}^2 = \abs{J_1}^2 + \eps^2\abs{J_2}^2 + 2\eps J_1\cdot J_2$, so we can use the inequality $2\abs{J_1\cdot J_2} \leq \abs{J_1}^2 + \abs{J_2}^2$
to obtain
\begin{align*}
    (1-\eps)( \abs{J_1}^2 + \eps\abs{J_2}^2)\leq \abs{\gradv H}^2 \leq (1+\eps)( \abs{J_1}^2 + \eps\abs{J_2}^2).
\end{align*}
In order to obtain the lower bound in \eqref{eq:gradH}, we drop the term $\abs{J_2}^2$ on the left hand side and we use $\abs{J_1}^2 = \ell^2 E^{2\ell-2}\abs{v}^2$. Concerning the upper bound in \eqref{eq:gradH}, we use
\begin{align*}
    \abs{J_2}^2 &=\abs*{\frac{\wangle{x}^{A}}{\wangle{v}^{B}} \, x - B \frac{\wangle{x}^{A} (x\cdot v) }{\wangle{v}^{B+2}} \, v}^2 \leq (1+B) \left( \frac{\wangle{x}^{A+1}}{\wangle{v}^{B}} \right)^2 \leq 2 \left(\wangle{x}^{\alpha(\ell-1)}\right)^2\leq 2c^2  E^{2(\ell-1)}
\end{align*}
so that,
\begin{align*}
    \left\vert \nabla_v H \right\vert^2 \leq(1+\eps)\big( \ell^2 E^{2\ell - 2} \wangle{v}^{2} + 2\eps c^2  E^{2(\ell-1)}\big)\lesssim E^{2\ell-2}\wangle{v}^{2}.
\end{align*}
This is exactly the upper bound we wanted to prove.

{\bf Proof of \eqref{eq:L*H}.} We just need to compute $\cL^* H$ term by term and perform some estimates. Concerning the energy term, we recall that $\sfT E =0$, so
\begin{align*}
    \cL^*(E^\ell) &= \Delta (E^\ell) + \frac{\gradv \M}{\M} \cdot \gradv (E^\ell) =\left[ (\ell-1) \frac{\abs{v}^2}{E} + d + v\cdot\frac{\gradv \M}{\M}\right]\ell E^{\ell-1}.
\end{align*}
Concerning the mixed term in $H$, we first notice that
\begin{align*}
    v \cdot \gradx\left(\frac{\wangle{x}^A}{\wangle{v}^B}x\cdot v\right)& = v \cdot \left( \frac{\wangle{x}^A}{\wangle{v}^B}v +A \frac{\wangle{x}^{A-2}(x\cdot v)x}{\wangle{v}^B} \right) = \frac{\wangle{x}^A}{\wangle{v}^B}  \left( \vert v \vert^2  + A \frac{(x\cdot v)^2}{\wangle{x}^{2}} \right)
\end{align*}
and, interchanging the roles of $A$ and $B$, we get similarly, 
\begin{align*}
    x \cdot \gradv\left(\frac{\wangle{x}^A}{\wangle{v}^B}x\cdot v\right) &= x \cdot \left( \frac{\wangle{x}^A}{\wangle{v}^B}x -B \frac{\wangle{x}^A(x\cdot v) \, v}{\wangle{v}^{B+2}}\right) = \frac{\wangle{x}^A}{\wangle{v}^B}  \left( \vert x \vert^2  - B \frac{(x\cdot v)^2}{\wangle{v}^{2}} \right).
\end{align*}
Next, we have
\begin{align*}
    \Delta_v\left(\frac{\wangle{x}^A}{\wangle{v}^B}x\cdot v\right) &=  -B\frac{\wangle{x}^A}{\wangle{v}^{B}}(x\cdot v)\left[  
    \frac{d+2}{\wangle{v}^{2}} -(B+2)\frac{\abs{v}^2}{\wangle{v}^{4}}\right].
\end{align*}
Gathering all terms, we can write 
\begin{align*}
&\frac{\wangle{v}^B}{\wangle{x}^A}\cL^*\left( \frac{\wangle{x}^A}{\wangle{v}^B} (x\cdot v) \right)\\
&= \vert v \vert^2  + A \frac{(x\cdot v)^2}{\wangle{x}^{2}} - \wangle{x}^{\alpha-2} \left( \vert x \vert^2  - B \frac{(x\cdot v)^2}{\wangle{v}^{2}}\right) \\
& \qquad \qquad \qquad \qquad -B (x\cdot v)\left[  
    \frac{d+2}{\wangle{v}^{2}} -(B+2)\frac{\abs{v}^2}{\wangle{v}^{4}}\right] + \left( x  -B \frac{(x\cdot v)}{\wangle{v}^{2}}v\right)\cdot\frac{\gradv\M}{\M}.
\end{align*}
We bound each term by the Cauchy-Schwarz inequality
\[
\abs{v}^2+A\frac{(x\cdot v)^2}{\wangle{x}^2}\leq (1+\abs{A})\abs{v}^2 \lesssim \abs{v}^2,
\]
\begin{equation*}
- \wangle{x}^{\alpha-2} \left( \vert x \vert^2  - B \frac{(x\cdot v)^2}{\wangle{v}^{2}}\right) \leq - (1-B) \, \wangle{x}^{\alpha-2}\vert x \vert^2,
\end{equation*}
\[
- B (x\cdot v)\left[ \frac{d+2}{\wangle{v}^{2}} -(B+2)\frac{\abs{v}^2}{\wangle{v}^{4}}\right] \leq \abs{B}(\abs{B}+d+4)\abs{x} \lesssim \abs{x},
\]
\begin{equation*}
 \left( x  -B \frac{(x\cdot v)}{\wangle{v}^{2}}v\right)\cdot\frac{\gradv\M}{\M} = (1+\abs{B})\abs{x}\abs*{\frac{\gradv\M}{\M}} \lesssim \abs{x} \wangle{v}^{\beta-1}.
\end{equation*}
In the last inequality includes both cases when $\M$ is given by \eqref{eq:Mexp} and \eqref{eq:Mpoly}, as for the latter case we have set $\beta=0$ by definition. Therefore, we can find a constant $C>0$ such that 
\begin{align*}
   \frac{\wangle{v}^B}{\wangle{x}^A}\cL^*\left( \frac{\wangle{x}^A}{\wangle{v}^B} (x\cdot v) \right) + (1-B)\abs{x}^2\wangle{x}^{\alpha-2} &\leq C \left[\abs{v}^2+ \wangle{v}^{(\beta -1)_+}  \abs{x}\right],
\end{align*}
which allows to conclude \eqref{eq:L*H}.

\noindent\textbf{\# Proof of \eqref{eq:L*Phi(H)}}. Finally, for any $\Phi\in\cC^2(\R_+)$, we have 
\begin{align*}
    \gradx (\Phi(H))&= \Phi'(H)\gradx H & \gradv (\Phi(H))&= \Phi'(H)\gradv H 
\end{align*}
and
\[
\Delta_v(\Phi(H))= \Phi''(H)\abs{\gradv H}^2 + \Phi'(H)\Delta_vH,
\]
so \eqref{eq:L*Phi(H)} follows.

\end{proof}

\subsection{Sub-Gaussian local equilibrium.}

The purpose of this section is to prove the following proposition, that provides a Lyapunov function with sub-exponential growth at infinity. 

\begin{prop}\label{prop:lyapEXP}
Assume that the external potential $V$ is given by \eqref{eq:V} with $\alpha>1$ and the local equilibrium is given by \eqref{eq:Mexp} with $\beta>0$. Consider the function $H$ defined as in \eqref{eq:H} with $\ell=2$, $A=a\alpha$, $B=1-b$
where $a,b\in(0,1)$. If $a,b,\eps > 0$ are small enough, then $H \asymp E^2$ and 
\begin{equation}\label{eq:LyapH}
    \mathcal{L}^*H\lesssim C\1_{\abs{v}\leq V_0,\,\abs{x}\leq X_0}-  H^{\frac{1}{2}}\wangle{v}^\beta
\end{equation}
for some constants $C,\,V_0,\,X_0>0$ large enough. As a consequence, for every $\theta\in(0,\min(1,\tfrac{\beta}{2})]$, there exists $\delta>0$ small enough such that the function $m=\exp(\delta H^{\frac{\theta}{2}})$ satisfies the Lyapunov condition
    \begin{equation}\label{eq:Lyapmsubexp}
        \cL^* m\lesssim C\1_{\abs{v}\leq V_0,\,\abs{x}\leq X_0} - m\,(\ln m)^{-\tfrac{1-\theta}{\theta}}.
    \end{equation}
for some constants $C,\,V_0,\,X_0>0$ large enough.
\end{prop}

\begin{proof}[{\bf Proof of \Cref{prop:lyapEXP}}]
Observe importantly that, in this setting, the conditions \eqref{eq:conAB} become
\begin{align*}
    0<&a\alpha<\alpha-1, & 0<b<1,
\end{align*}
which are satisfied when $a$ and $b$ are small enough. Therefore, by \Cref{lem:tech}, we have
\begin{align*}
    H&\asymp E^2 & \abs{\gradv H}^2\lesssim E^{2}\wangle{v}^2.
\end{align*}

\medskip 

\noindent\textbf{\# Proof of \eqref{eq:LyapH}}. Thanks to \eqref{eq:L*H}, substituting $\tfrac{\gradv\M}{\M}=-\wangle{v}^{\beta-2}v$ we find
\begin{align*}
   \cL^*H &\leq 2\left[ \frac{\abs{v}^2}{E} + d - \abs{v}^2\wangle{v}^{\beta-2} \right]  E - b\eps\frac{\wangle{x}^{\alpha a + \alpha-2} \abs{x}^2}{\wangle{v}^{1-b}} + C\eps \frac{\wangle{x}^{\alpha a}}{\wangle{v}^{1-b}} \left[ \abs{v}^2 +\wangle{v}^{(\beta-1)_+}  \abs{x} \right] \\
    &\leq  2\left[ (2+d)\1_{\vert v \vert \leq V_0} -  \frac12 \wangle{v}^{\beta} \right] E - b\eps\frac{\wangle{x}^{\alpha a + \alpha-2} \abs{x}^2}{\wangle{v}^{1-b}} + C\eps \wangle{v}^{b-1}\wangle{x}^{a\alpha} \left[ \wangle{v}^{2} + \wangle{v}^{(\beta-1)_+}\wangle{x} \right],
\end{align*}
where we used $\frac{\abs{v}^2}{E}\leq 2$ and we took $V_0>0$ large enough. We bound the last term as follows
\begin{align*}
  \wangle{v}^{b-1}\wangle{x}^{a\alpha} \left[ \wangle{v}^{2} + \wangle{v}^{(\beta-1)_+}\wangle{x} \right]&\lesssim\wangle{v}^{b-1}E^{a} \left[ \wangle{v}^{2} + \wangle{v}^{(\beta-1)_+}E^{\frac{1}{\alpha}} \right] \\
  &=  \wangle{v}^{\beta} E \left[  \wangle{v}^{b  +1- \beta}E^{a-1} +  \wangle{v}^{b  -1 +(\beta-1)_+ - \beta}E^{a +\frac{1}{\alpha}-1}  \right]
\end{align*}
Notice that $(\beta-1)_+ -\beta\leq 0$ for all $\beta\geq 0$ and, if $a,b>0$ are small enough, we also have $b +2a \leq \beta +1$ and $a\leq 1-\frac{1}{\alpha}$, so we conclude
\begin{align*}
    \wangle{v}^{b-1}\wangle{x}^{a\alpha} \left[ \wangle{v}^{2} + \wangle{v}^{(\beta-1)_+}\wangle{x} \right]\lesssim \wangle{v}^{\beta} E \left[  \wangle{v}^{b - \beta+2a-1} +  \wangle{v}^{b - 1}E^{a +\frac{1}{\alpha}-1}  \right] \lesssim \wangle{v}^{\beta} E.
\end{align*}
Moreover,
\begin{align*}
2(2+d) \1_{\vert v \vert \leq V_0} E - b\eps\frac{\wangle{x}^{\alpha a + \alpha-2} \abs{x}^2}{\wangle{v}^{1-b}} &\leq 2(2+d) \1_{\vert v \vert \leq V_0} E  - b\eps  \1_{\vert v \vert \leq V_0}\frac{\wangle{x}^{\alpha a + \alpha-2} \abs{x}^2}{\wangle{v}^{1-b}}\\
&\leq  \1_{\vert v \vert \leq V_0} \wangle{v}^{b-1} \left( 2(2+d) \wangle{v}^{1-b}E - b\eps\wangle{x}^{\alpha a + \alpha-2} \abs{x}^2 \right)\\
& \lesssim \1_{\vert v \vert \leq V_0}\1_{\vert x \vert \leq X_0} ,
\end{align*}
when $X_0$ is large enough because $\wangle{v}^{1-b}E\lesssim \wangle{x}^\alpha$ if $\abs{v}\leq V_0$. Collecting all the estimates we have that, for some constant $C_1>0$,
\[
\cL^* H\leq C_1 \1_{\vert v \vert \leq V_0}\1_{\vert x \vert \leq X_0} -\wangle{v}^\beta E + C_1\eps \wangle{v}^\beta E
\]
which gives the conclusion  after taking $\eps>0$ small enough and rememberin that $E\asymp H^{\frac12}$.

\medskip\noindent\textbf{\# Proof of \eqref{eq:Lyapmsubexp}}. Let us now consider the function $m= \Phi\left( H \right) := \exp(\delta H^{\frac{\theta}{2}})$,  $\theta\in \left(0,\min\left(1,\tfrac{\beta}{2}\right)\right]$ and $\delta>0$. Recalling \Cref{lem:tech}, we have
\begin{equation*}
   \mathcal{L}^*m = \Phi'(H) \mathcal{L}^*(H) + \Phi''(H)\abs{\gradv H}^2.
\end{equation*}
After direct computations, 
\begin{align*}
    &\Phi'(H) = \frac{\delta \theta}{2} H^{\frac{\theta}{2} -1} \Phi(H) = \frac{\delta \theta}{2} H^{\frac{\theta}{2} -1} m\\
    &\Phi''(H) = \frac{\delta \theta}{2} H^{\frac{\theta}{2} -1} \Phi'(H)   + \frac{\delta \theta}{2} \left(\frac{\theta}{2} -1\right) H^{\frac{\theta}{2} -2}\Phi(H) \leq \delta H^{\frac{\theta}{2} -1} \Phi'(H)
\end{align*}
the latter since $\frac{\theta}{2} -1 \leq 0$. We end up with, 
\begin{align*}
    \mathcal{L}^*m&= \Phi'(H) \mathcal{L}^*(H) + \Phi''(H)\abs{\gradv H}^2 \leq \Phi'(H) \left( \mathcal{L}^*(H) + \delta H^{\frac{\theta}{2} -1}\abs{\gradv H}^2\right).
\end{align*}
Thanks to \eqref{eq:gradH}, we have
\[
H^{\frac\theta2 -1}\abs{\gradv H}^2 \lesssim H^{\frac\theta2 -1} E^{2}\wangle{v}^2\lesssim H^{\frac\theta2 }\wangle{v}^2,
\]
hence, for some suitable constant $C_2$
\begin{align*}
    \mathcal{L}^*m&\lesssim \Phi'(H) \left( C\1_{\abs{v}\leq V_0,\,\abs{x}\leq X_0} - H^{\frac{1}{2}}\wangle{v}^\beta + C_2\delta H^{\frac{\theta}{2} }\wangle{v}^2\right)\\
    &\leq \Phi'(H) \left( C\1_{\abs{v}\leq V_0,\,\abs{x}\leq X_0} - H^{\frac{1}{2}}\wangle{v}^\beta \left( 1 - C_2\delta  H^{\frac{\theta}{2} -\frac{1}{2}} \wangle{v}^{2-\beta} \right) \right).
\end{align*}
If $\beta\geq 2$, then we immediately have $H^{\frac{\theta}{2} -\frac{1}{2}} \wangle{v}^{2-\beta}\lesssim1$ because $\theta\leq 1$. On the other hand, if $\beta\in(0,2)$, then we have 
\[
H^{\frac{\theta}{2} -\frac{1}{2}} \wangle{v}^{2-\beta}\lesssim H^{\frac{\theta}{2} -\frac{1}{2}} E^{1-\frac{\beta}{2}}\lesssim H^{\frac{\theta}{2}- \frac\beta4}\lesssim1, 
\]
because $\theta\leq \frac{\beta}{2}$. In any case, by taking $\delta>0$ small enough, we have
\begin{align*}
    \cL^*m\lesssim \Phi'(H) \left( C\1_{\abs{v}\leq V_0,\,\abs{x}\leq X_0} - H^{\frac{1}{2}}\wangle{v}^\beta \right).
\end{align*}
Finally, coming back to $m$,
\begin{align*}
    \mathcal{L}^*m&\lesssim \Phi'(H) \left( C\1_{\abs{v}\leq V_0,\,\abs{x}\leq X_0} - H^{\frac{1}{2}}\wangle{v}^\beta \right) \lesssim  H^{\frac{\theta}{2} -1} m\left( C\1_{\abs{v}\leq V_0,\,\abs{x}\leq X_0} - H^{\frac{1}{2}}\wangle{v}^\beta \right)  \\
    & \lesssim  C\1_{\abs{v}\leq V_0,\,\abs{x}\leq X_0}-H^{\frac{\theta - 1}{2}}\wangle{v}^\beta m\lesssim C\1_{\abs{v}\leq V_0,\,\abs{x}\leq X_0}- m (H^\frac{\theta}{2})^{\frac{\theta - 1}{\theta}}\\
    &\lesssim C\1_{\abs{v}\leq V_0,\,\abs{x}\leq X_0}- m\,(\ln m)^{-\frac{1-\theta}{\theta}},
\end{align*}
which concludes the proof.
\end{proof}

\subsection{Fat-tailed local equilibrium}

\begin{prop}\label{prop:lyapEXPpoly}
Assume that the external potential $V$ is given by \eqref{eq:V} for some $\alpha>1$ and the local equilibrium is given by \eqref{eq:Mpoly} with $\gamma>0$. Consider the function $H$ defined as in \eqref{eq:H} with $\ell\geq \max\left(2,1+\frac\gamma2\right)$, $A=\alpha(\ell-2 +a)$, $B=1-b$
where $a,b\in(0,1)$. If $a,b,\eps > 0$ are small enough, then 
\begin{align}\label{eq:boundHpoly}
    H&\asymp E^\ell & (1-\eps) \ell^2 E^{2\ell-2}\abs{v}^2\leq &\left\vert \nabla_v H \right\vert^2
\end{align}
and, for every $k\in\left(1,1+\frac\gamma2 \right)$, the function $m=H^{\frac{k}{\ell}}$ satisfies the Lyapunov condition
    \begin{align}\label{eq:LyapHpoly}
        \cL^* m\lesssim C\1_{\abs{v}\leq V_0,\,\abs{x}\leq X_0}  - m^{1-\tfrac{1}{k}},
    \end{align}
for some constants $C,\,V_0,\,X_0>0$ large enough.
\end{prop}

\begin{proof}[{\bf Proof of \Cref{prop:lyapEXPpoly}}]

First of all, we notice that \eqref{eq:conAB} are satisfied for our choice of $A$ and $B$, if $a$ and $b$ are small enough. Therefore \Cref{lem:tech} ensures that
\begin{align*}
    (1- c \eps)E^\ell \leq& H\leq (1+ c \eps)E^\ell &
    (1-\eps) \ell^2 E^{2\ell-2}\abs{v}^2\leq &\left\vert \nabla_v H \right\vert^2 
\end{align*}
for $\eps$ sufficiently small.

\medskip\noindent\textbf{\# Proof of \eqref{eq:LyapHpoly}}.
Let $m=H^{\tfrac{k}{\ell}}$ with $k\in\left(1,1+\frac\gamma2 \right)$; notice that $m$ is a concanve function of $H$ because $k\leq \ell$. Then, thanks to \eqref{eq:L*Phi(H)}, \eqref{eq:boundH} and \eqref{eq:gradH}, one has,
\begin{align*}
    \cL^* m &= \frac{k}{\ell}H^{\frac{k}{\ell}-1}\Big(\cL^* H - \Big(1-\tfrac k \ell \Big) \frac{\abs{\gradv H}^2}{H}\Big)\leq \frac{k}{\ell}H^{\frac{k}{\ell}-1}\Big(\cL^* H - \frac{(\ell- k)(1-\eps)\ell }{1+c\eps} E^{\ell-2}\abs{v}^2\Big) \\
\end{align*}
We compute $\cL^*H$ by plugging $\tfrac{\gradv\M}{\M}=-(d+\gamma)\wangle{v}^{-2}v$ into \eqref{eq:L*H}, which gives us
\begin{align*}
    \cL^*H&\leq \left[ (\ell-1)\frac{\abs{v}^2}{E} + d - (d+\gamma)\abs{v}^2\wangle{v}^{-2} \right] \ell E^{\ell-1} - b\eps\frac{\wangle{x}^{\alpha(\ell -1+a) -2}\abs{x}^2 }{\wangle{v}^{1-b}}  + C\eps\frac{\wangle{x}^{\alpha(\ell -2+a)} }{\wangle{v}^{1-b}} \left(  \abs{v}^2 + \abs{x} \right) .
\end{align*}
As a consequence, after reorganizing the terms, we obtain the following estimate
\begin{align*}
    \cL^*m&\leq \frac{k}{\ell}H^{\frac{k}{\ell}-1}\Big(\left[ \left(\ell - \frac{(\ell- k)(1-\eps) }{1+c\eps}-1\right)\frac{\abs{v}^2}{E} + d - (d+\gamma)\abs{v}^2\wangle{v}^{-2} \right] \ell E^{\ell-1} \\
    &\qquad\qquad\qquad - b\eps\frac{\wangle{x}^{\alpha(\ell -1+a) -2}\abs{x}^2 }{\wangle{v}^{1-b}}  + C\eps\frac{\wangle{x}^{\alpha(\ell -2+a)} }{\wangle{v}^{1-b}} \left(  \abs{v}^2 + \abs{x} \right)\Big).
\end{align*}
Now we need to simplifying this expression. First, we notice that
\begin{align*}
    \ell - \frac{(\ell- k)(1-\eps) }{1+c\eps}-1 &= \ell\Big(1-\frac{1-\eps}{1+c\eps}\Big) + k\frac{1-\eps}{1+c\eps} -1\\
    &\leq \frac{(1+c)\ell\eps}{1+c\eps}  +k -1\leq (1+c)\ell\eps +k -1 .
\end{align*}
Combining this estimate with $\frac{\abs{v}^2}{E}\leq 2$ and $\abs{v}^2\wangle{v}^{-2} = 1 - \wangle{v}^{-2}$ we get 
\[
\left(\ell - \frac{(\ell- k)(1-\eps) }{1+c\eps}-1\right)\frac{\abs{v}^2}{E} + d - (d+\gamma)\abs{v}^2\wangle{v}^{-2} \leq 2(1+c)\ell\eps +2(k -1) -\gamma + (d+\gamma)\wangle{v}^{-2} 
\]
Since $k< 1+\tfrac{\gamma}{2}$, we can take $\eps>0$ small enough so that
\[
\eta:=\gamma- 2(k - 1) - 2(1+c)\ell\eps>0.
\]
Hence we obtain 
\begin{align*}
    \cL^*m&\leq \frac{k}{\ell}H^{\frac{k}{\ell}-1}\Big( -\eta\ell E^{\ell-1} + (d+\gamma)\ell\wangle{v}^{-2} E^{\ell-1}- b\eps\frac{\wangle{x}^{\alpha(\ell -1+a) -2}\abs{x}^2 }{\wangle{v}^{1-b}}  + C\eps\frac{\wangle{x}^{\alpha(\ell -2+a)} }{\wangle{v}^{1-b}} \left(  \abs{v}^2 + \abs{x} \right)\Big).
\end{align*}
The next step is to bound the error terms. Observe that
\begin{equation*}
\frac{\wangle{x}^{\alpha(\ell-2+a)} }{\wangle{v}^{1-b}}   \abs{v}^2  \leq  \wangle{x}^{\alpha(\ell-2+a)} \wangle{v}^{1+b} \lesssim E^{\ell-2 + a }E^{ \frac{1+b}{2}} \leq E^{\ell-1},
\end{equation*}
provided $a + \frac{1+b}{2} \leq 1$, which is true for $a,\,b>0$ small enough. Moreover, 
\begin{equation*}
   \frac{\wangle{x}^{\alpha(\ell-2+a)} }{\wangle{v}^{1-b}} \abs{x}  \leq \wangle{x}^{\alpha(\ell-2)+\alpha a +1} \lesssim E^{\ell-2 + \frac{\alpha a+1}{\alpha}} \leq  E^{\ell-1},
\end{equation*}
when $a + \frac{1}{\alpha} \leq 1$. Consequently, by taking $\eps,a,b>0$ small and $V_0,X_0,C_3>0$ large enough, we have
\begin{align*}
    \cL^*m    &\lesssim H^{\frac{k}{\ell}-1}\Big( -\frac{\eta\ell}{2} E^{\ell-1} + (d+\gamma)\ell\wangle{v}^{-2} E^{\ell-1}- b\eps\frac{\wangle{x}^{\alpha(\ell -1+a) -2}\abs{x}^2 }{\wangle{v}^{1-b}}  \Big)\\
    &\leq H^{\frac{k}{\ell}-1}\Big(\left(  - \frac{\eta\ell}{4} + (d+\gamma)\ell \1_{\vert v \vert \leq V_0} \right) E^{\ell-1} -b\eps \1_{\vert v \vert \leq V_0}\frac{\wangle{x}^{\alpha(\ell-1+a)-2}\abs{x}^2 }{\wangle{v}^{1-b}}\Big)\\
    &\leq H^{\frac{k}{\ell}-1}\Big( -\frac{\eta\ell}{4}E^{\ell-1} +  \1_{\vert v \vert \leq V_0}\left( (d+\gamma)\ell E^{\ell-1}  -b\eps\frac{\wangle{x}^{\alpha(\ell-1+a)-2}\abs{x}^2 }{\wangle{v}^{1-b}} \right)  \Big)\\
    &\leq H^{\frac{k}{\ell}-1}\Big( -\frac{\eta\ell}{4}E^{\ell-1} +  C_3\1_{\vert v \vert \leq V_0}\1_{\vert x \vert \leq X_0}  \Big),
\end{align*}
because $\wangle{v}^{1-b}E^{\ell-1}\lesssim \wangle{x}^{\alpha(\ell-1)}$ if $\abs{v}\leq V_0$. We conclude that
\begin{align*}
    \cL^*m \lesssim H^{\frac{k}{\ell}-1}\big( C_3\1_{\vert v \vert \leq V_0}\1_{\vert x \vert \leq X } - E^{\ell-1}\big)\lesssim  C_3\1_{\abs{v}\leq V_0,\,\abs{x}\leq X_0} - H^{\frac{k-1}{\ell}} \lesssim C_3\1_{\abs{v}\leq V_0,\,\abs{x}\leq X_0} - m^{\frac{k-1}{k}}
\end{align*}
that is exactly \eqref{eq:LyapHpoly}.

\end{proof}

\section{About the positivity condition}\label{sec:positivity}
In this section, we prove the positivity condition needed to apply the Harris theorem. We follow the idea of \cite{CGMM24}, which gets the positivity condition as a direct consequence of the Harnack inequality. We denote by $S_\cL$  the semigroup associated to \eqref{L}.

\begin{prop}\label{prop:Harnack}
For any $0<T_0<T_1$ and compact set $\mathcal{C}\subseteq \RRd$, there holds
\begin{equation}\label{eq:Harnack}
    \sup_{\mathcal{C}} S_\cL(T_0)f_0 \leq C\inf_{\mathcal{C}} S_\cL(T_1)f_0,
\end{equation}
for some constant $C = C(T_0, T_1, \mathcal{C}) > 0$.
\end{prop}
The proof is based on a local version of the Harnack inequality and an iteration argument to obtain \eqref{eq:Harnack}. The local Harnack inequality has been obtained in \cite{GM22} by De Giorgi methods, see also \cite{MR3923847} for previous results, \cite[Theorem~5.1]{bony_1969} for seminal considerations in this topic, together with the Harnack inequality proved by Lanconelli and Polidoro in \cite[Theorem~5.1]{Lanconelli_1994}. The iteration argument is taken from \cite{LRR22}, and is an adaptation of the ideas from \cite{MR4017782}. We recall the main result of \cite{GM22}. Consider the following partial differential equation
\begin{equation}\label{eq:GueMou}
    \partial_t g =\mathfrak{M}g :=- v \cdot \gradx g + \Delta_v g + b \cdot \gradv g ,\qquad  t \in \R, x \in \R^d, v \in \R^d,
\end{equation}
where $b=b(x,v)$ is a measurable and locally bounded vector field. The class of equations \eqref{eq:GueMou} is invariant under the set of transformations 
\[
\cT_{\eta,z_0}\colon(t, x, v) \mapsto (t_0 + \eta^2t, x_0 + \eta^3x + \eta^2tv_0, v_0 + \eta v),
\]
where $r>0$ and $z=(t_0, x_0, v_0)\in \R\times\RRd$. For this reason, we define for $z_0 \in \R\times\RRd$ and $r > 0$ the set
\[
Q_r(z_0) := \Big\{-r^2 < t - t_0 \leq 0, \abs{x - x_0 - (t - t_0)v_0} < r^3
, \abs{v - v_0} < r	\Big\}
\] %z_0 \circ [rQ_1] :=
and observe that $\cT_{\eta,z_0}(Q_r(0,0,0))= Q_{\eta r}(z_0)$.  The following theorem is taken from \cite[Thm. 5 \& Prop. 12]{GM22} but it is rewritten as in \cite{CGMM24} using the invariant transformations and the local boundedness of $b$.

\begin{thm}\label{thm:GueMou}
Let $z_0=(t_0,x_0,v_0)$ with $t_0\in\R$ and $(x_0,v_0)\in \cC\subseteq\RRd$ a compact subset. We have the following two results.
\begin{enumerate}[label=(\roman*)]
    \item There is $\zeta > 0$ such that any non-negative weak super-solution $g$ to \eqref{eq:GueMou} in $Q_1(z_0)$ satisfies %the weak Harnack inequality
    \[
    \norm{g}_{\sfL^\zeta(\widetilde{Q}_{r_0}^-(z_0))} \lesssim \inf_{Q_{r_0}(z_0)}g
    \]
    where $r_0 = 1/40$ and $\widetilde{Q}^-_{r_0}(z_0):=\cT_{1,z_0}(Q_{r_0}(-\tau,0,0)) =  Q_{r_0}(z_0-(\tau, \tau v_0, 0))$ with $\tau = \frac{19}{2}r^2_0$.
    \item Let $h$ be a non-negative weak sub-solution to \eqref{eq:GueMou} in an open set $\mathcal{U}\subseteq \R^{1+2d}$. Given any $Q_r(z_0) \subseteq Q_R(z_0) \subseteq\mathcal{U}$ with $0 < r < R\leq 1$, and $\zeta> 0$, $h$ satisfies
    \begin{equation}
    \norm{h}_{\sfL^\infty(Q_r(z_0))}\lesssim C \norm{h}_{\sfL^\zeta(Q_R(z_0))}
    \end{equation}
    with $C=C(R,r,v_0, \zeta)$.
\end{enumerate}

\end{thm}
\begin{proof}[{\bf Proof of \Cref{prop:Harnack}}]
\textit{\# Step 1: Local Harnack inequality}.
We prove that for any $f$ solution to \eqref{L}, i.e. such that
\[
\dt f = -v\cdot\gradx f + \gradx V\cdot\gradv f +\Delta_v f - \frac{\gradv\M}{\M}\cdot \gradv f - \gradv\cdot \left(\frac{\gradv\M}{\M}\right)f,
\]
and for any $z_0\in\R^{1+2d}$ there holds
\begin{equation}\label{eq:localHarnack}
    \sup_{\widetilde{Q}^-_{r_0/2}(z_0)}f\lesssim \inf_{Q_{r_0}(z_0)}f,
\end{equation}
where $r_0=1/40$ and $\widetilde{Q}^-_{r_0/2}(z_0)= Q_{r_0/2}(z_0 - (\tau,\tau v_0,0))$ with $\tau = \frac{19}{2}r_0^2$. Notice that for $\M$ as in \eqref{eq:Mexp} or as in \eqref{eq:Mpoly} we have
\[
\norm*{\gradv\cdot \left(\frac{\gradv\M}{\M}\right)}_{\sfL^\infty(\cC)}=:\lambda<\infty .
\]
We define $g:=e^{\lambda t}f$ and we consider \eqref{eq:GueMou} with $b(x,v)=\gradx V(x) - \tfrac{\gradv \M(v)}{\M(v)}$. We have
\begin{align*}
    \dt g = \mathfrak{M}g + \left(\lambda - \gradv\cdot \left(\frac{\gradv\M}{\M}\right)\right)g\geq  \mathfrak{M}g,
\end{align*}
in $\R^{1+2d}$ so that $g$ is a non-negative weak super-solution to \eqref{eq:GueMou}. From \Cref{thm:GueMou}-\textit{(i)} we deduce that
\[
\norm{g}_{\sfL^\zeta(\widetilde{Q}_{r_0}^-(z_0))}\lesssim \inf_{Q_{r_0}(z_0)}g
\]
which, minimizing $t\mapsto e^{\lambda t}$ over $\widetilde{Q}_{r_0}^-(z_0)$ and maximizing it over $Q_{r_0}(z_0)$, rewrites as
\begin{equation}\label{Har1}
\norm{f}_{\sfL^\zeta(\widetilde{Q}^-_{r_0}(z_0))}\lesssim\inf_{Q_{r_0}(z_0)}f.
\end{equation}

On the other hand, let us define $h=e^{-\lambda t}f$ and we notice that 
\begin{align*}
    \dt h = \mathfrak{M}h - \left(\lambda + \gradv\cdot \left(\frac{\gradv\M}{\M}\right)\right)h\leq  \mathfrak{M}h \qquad \text{in } \R^{1+2d}
\end{align*}
Thus $h$ is a non-negative sub-solution of \eqref{eq:GueMou} and by applying \Cref{thm:GueMou}-\textit{(ii)} with $r=r_0/2$, $R=r_0$ and with $z_0-(\tau, \tau v_0, 0)$ in place of $z_0$ we deduce
\begin{equation*}
    \norm{h}_{\sfL^\infty(\widetilde{Q}^-_{r_0/2}(z_0))}\lesssim  \norm{h}_{\sfL^\zeta(\widetilde{Q}_{r_0}^-(z_0))}.
\end{equation*}
By minimizing $t\mapsto e^{-\lambda t}$ over $\widetilde{Q}^-_{r_0/2}(z_0)$ and maximizing it over $\widetilde{Q}_{r_0}^-(z_0)$ we obtain 
\begin{equation*}
    \norm{f}_{\sfL^\infty(\widetilde{Q}^-_{r_0/2}(z_0))}\lesssim  \norm{f}_{\sfL^\zeta(\widetilde{Q}_{r_0}^-(z_0))},
\end{equation*}
which with \eqref{Har1} gives \eqref{eq:localHarnack}.

\textit{\# Step 2: Proof of \eqref{eq:Harnack}}. Once the local Harnack inequality \eqref{eq:localHarnack} is satisfied, we can deduce \eqref{eq:Harnack} following the second step in the proof of \cite[Thm 2.15]{LRR22}. The main idea is to repeatedly apply  \eqref{eq:localHarnack} along a chain of sets $Q_r(z_k)$ with $r>0$ sufficiently small and $(z_k)_{k\geq0}\subseteq\R\times\RRd$. The points $(z_k)_{k\geq0}$ can be constructed explicitly with the property  $z(s_k)\in\widetilde{Q}^-_{r}(z(s_{k+1}))$. In this way, for any $(x,v),(x',v')\in \cC$ and $f=S_\cL(t)f_0$ we can find a constant $C$ such that
\[
f(T_0,x,v)\leq Cf(T_1,x,v)
\]
where $C$ only depends on $T_0$, $T_1$ and $\cC$.

\end{proof}

We are now ready to prove the positivity condition we need for the Harris theorem. For any $M>0$, we consider the compact sets
\begin{equation*}
    \cC =\{(x,v)\in\RRd \,:\, m(x,v)\leq M\}
\end{equation*}
where $m$ is a Lyapunov function given by \Cref{prop:lyapEXP} or \Cref{prop:lyapEXPpoly}. In particular, notice that $\cC$ is compact for any $M>0$.

\begin{lemma}\label{lem:positivity}
For any $M > 0$ and for any $0<T_0<T_1$, there exists a nonzero measure $\nu$ such that
\[
S_\mathcal{L}(T_1) f_0 \geq \nu\int_{\mathcal{C}}S_\cL(T_0) f_0,
\]
for all $f_0 \in \sfL^1(m)$, $f_0 \geq 0$.
\end{lemma}
\begin{proof}[{\bf Proof of \Cref{lem:positivity}}]
We know that for any $0<T_0<T_1$ and any compact set $\cC$ there exists a constant $C>0$ such that 
\[
\sup_\cC S_\cL(T_0)f_0\leq C\inf_\cC S_\cL(T_1)f_0
\]
We deduce that
\begin{align*}
    S_\cL(T_1)f_0 \geq \1_\cC \inf_\cC S_\cL(T_1)f_0 \geq \frac{\1_\cC}{C}(\sup_\cC S_\cL(T_0)f_0) \geq \frac{\1_\cC}{C\abs{\cC}} \int_\cC S_\cL(T_0)f_0.
\end{align*}
So the conclusion follows by taking $\nu= \frac{\1_\cC}{C\abs{\cC}}$.
\end{proof}

\section{Proofs of \Cref{thm:mainEXP} and \Cref{thm:mainPOLY}}\label{sec:proofs}

The positivity condition for the Harris theorem has been proved in \Cref{lem:positivity} for all cases, and the Lyapunov functions have been found in \Cref{prop:lyapEXP} and \Cref{prop:lyapEXPpoly}. 

\medskip\noindent\textbf{\# Proof of \Cref{thm:mainEXP}}.
Consider the function $m=\exp(\delta H^{\frac{\theta}{2}})$ defined as in \Cref{prop:lyapEXP}, then we have the Lyapunov condition
\[
\cL^*m\lesssim C\1_{B_R} - m (\ln m)^{-\frac{1-\theta}{\theta}}.
\]
The Harris theorem and the fact that $H\asymp E^2$ give existence and uniqueness of a steady state
\[
G\in\sfL^1 ( \exp(\delta E^{\frac{\beta}{2}}))
\]
and, for any $\theta\in(0,\min\big(\tfrac{\beta}{2},1\big)]$, the existence of some $\lambda>0$ such that
\begin{equation*}
    \left\Vert f(t)- \left(\int_{\R^d \times \R^d} f_0\right)G\right\Vert_{\sfL^1}\lesssim 
    e^{-\lambda t^\theta} \left\Vert f_0- \left(\int_{\R^d \times \R^d} f_0 \right) G\right\Vert_{\sfL^1(e^{\delta E^\theta})}.
\end{equation*}
for all $t\geq0$ and solutions $f$ with initial condition $f_0 \in \sfL^1(e^{\delta E^\theta})$.

\medskip\noindent\textbf{\# Proof of \Cref{thm:mainPOLY}}.
Consider now the function $m=H^{\frac{k}{\ell}}$ defined as in \Cref{prop:lyapEXPpoly}, then we have the Lyapunov condition 
\[
\cL^*m\lesssim C\1_{B_R}- m^{1-\frac{1}{k}}.
\]
The Harris theorem and the fact that $m\asymp E^k$ implies that there exists a steady state $G\in\sfL^1(E^{k-1})$ for all $k\in\left(1,1+\frac\gamma2\right)$ and that
\begin{equation*}
    \left\Vert f(t)- \left(\int_{\R^d \times \R^d} f_0 \right) G\right\Vert_{\sfL^1}\lesssim \frac{1}{\bangle{t}^{k}}\left\Vert f_0-\left(\int_{\R^d \times \R^d} f_0 \right)G\right\Vert_{\sfL^1(E^k)}
\end{equation*}
for all $t\geq 0$ and solutions $f$ with initial condition $f_0 \in \sfL^1(E^k)$.

\section{Numerical results}

This section is devoted to numerical simulations of the kinetic Fokker-Planck equation \eqref{L}. In particular, we are interested in describing the shape of the steady state for some choices of parameter $\alpha$ and $\beta$ in the strip $(\alpha,\beta)\in[1,\infty)\times(0,2)$.

\subsection{Description of the algorithm}
The simulations in this work are obtained using the finite volume numerical scheme presented in \cite{RB24} and the code is written in Python\footnote{The code is available at \url{https://github.com/luca-ziviani/KineticFokkerPlanck_Python}}. For simplicity, we restrict ourselves to the one-dimensional case $d=1$. We considered a rectangular domain $(x,v) \in [-L,L]\times[-V,V]$, with $L,V>0$ large enough, and we divided the intervals into cells with a uniform subdivision of size $\Delta x = \tfrac{2L}{N_x}$ and $\Delta v = \tfrac{2V}{N_v}$, where $N_x$ and $N_v$ are even positive numbers. The coordinates of the cell centres are $x_n= - L +(n+\tfrac{1}{2})\Delta x $ for $n=0,\dots, N_x-1$ and similarly $v_m= - V +(m+\tfrac{1}{2})\Delta v $ for $m=0,\dots, N_v-1$. The volume element is 
\[
\omega_{n,m} = [x_n-\tfrac{\Delta x}{2},x_n+\tfrac{\Delta x}{2}]\times[v_m-\tfrac{\Delta v}{2},v_m+\tfrac{\Delta v}{2}], 
\]
whose volume is $\abs{\omega_{n,m}} = \Delta x\Delta v$. Concerning the time coordinate, we fix $T>0$ large and partition the interval $[0,T]$ into $N_t>1$ sub-intervals of length $\Delta t = \tfrac{T}{N_t}$, so we have the points $t_k = k \Delta t$, for $k= 0,\dots,N_t$. We approximate a function $f=f(t,x,v)$ on points $(x_n,v_m)$ at time $t_k$ by the average
\[
f_{n,m}^k = \frac{1}{\abs{\omega_{n,m}}}\int_{\omega_{n,m}} f(t_k,x,v)\dx\dv.
\]
The main idea of a finite volume scheme is to compute the flux through each side of every element volume $\omega_{n,m}$ and use it to compute $\{f_{n,m}^{k+1}\}_{n,m}$ from the previous state $\{ f_{n,m}^{k}\}_{n,m}$. In \cite{RB24}, the author shows how to compute fluxes in such a way that the numerical scheme is explicit and stable under suitable CFL conditions. We briefly describe the scheme. For the advection-in-space term $v\cdot\gradx f$, the Kurganov-Tadmor flux \cite{KT00} is implemented, for the remaining differential operator in the velocity variable a generalised Chang-Cooper \cite{CC70} finite volume method is used. The boundary conditions are imposed on the fluxes through the exterior edges of the boundary volumes, and they are chosen to be specular on the sides $\{-L,L\}\times[-V,V]$ and zero-flux on the sides $[-L,L]\times\{-V,V\}$. Finally, the two types of flux are put together by a Strang splitting, and the evolution in time is performed by a Runge-Kutta method of order 2.

Since we are trying to simulate an unbounded domain $\R\times\R$ as a very large rectangle $[-L,L]\times[-V,V]$, we expect to see some perturbation of the solution for points near the boundary. Nevertheless, we hope to see our results confirmed at least for points relatively close to the origin.

\subsection{Plot of the density of the steady state}
In this section, we show some numerical results about the density $\rho_f$ of the solution $f$. The result of \Cref{thm:mainEXP} suggests that the steady state $G$ may resemble the function
\[
\widetilde{G}(x,v)= \exp(-\delta E^{\frac{\beta}{2}})
\]
for some $\delta>0$. Notice that the density $\rho_{\widetilde{G}}(x)=\int_\R \exp(-\delta E(x,v)^{\frac{\beta}{2}})\dv$ of such a function has the following asymptotic behaviour
\begin{equation}\label{eq:asympRhoG}
    \rho_{\widetilde{G}}(x) \sim C \abs{x}^{\tfrac{\alpha}{2}\left(1-\tfrac{\beta}{2}\right)}\exp\left(-\delta\left(\frac{\wangle{x}^\alpha}{\alpha}\right)^{\frac{\beta}{2}} \right) 
\end{equation}
as $\abs{x}\to \infty$, where $C= \frac{2\sqrt{\pi}\alpha^{\frac{\beta}{4}}}{\sqrt{\beta\delta\alpha}}$\footnote{This asymptotic result can be obtained as in \cite{BEZ25} by means of Laplace-type integrals, see for example \cite{T14}.}. Hence, we compare this asymptotic behaviour (with some $\delta>0$ to be determined) with the graph of $\rho_f$ at time $t\in[0,T]$ with $T$ large enough.

\begin{figure}
  \centering
  % Line 1
  \begin{subfigure}{0.48\textwidth}
    \includegraphics[width=\linewidth]{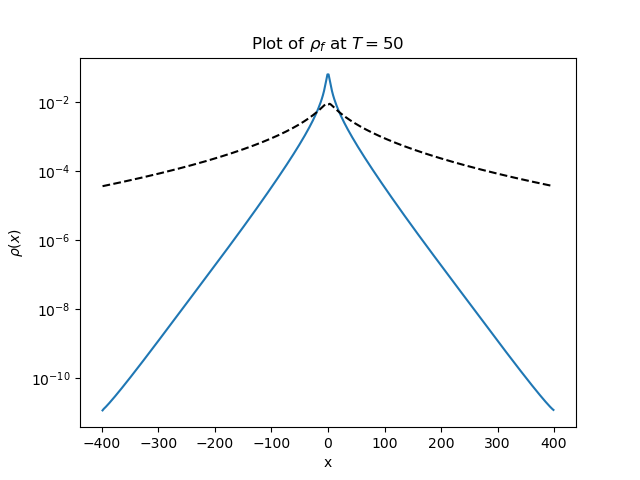}
  \end{subfigure}
  \begin{subfigure}{0.48\textwidth}
    \includegraphics[width=\linewidth]{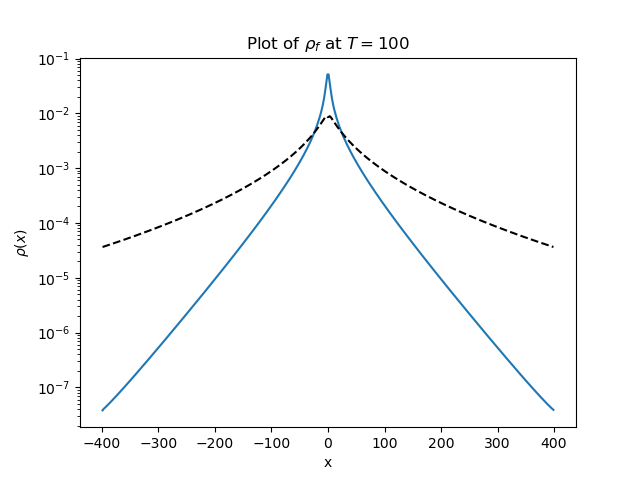}
  \end{subfigure}
  
  % Line 2
  \begin{subfigure}{0.48\textwidth}
    \includegraphics[width=\linewidth]{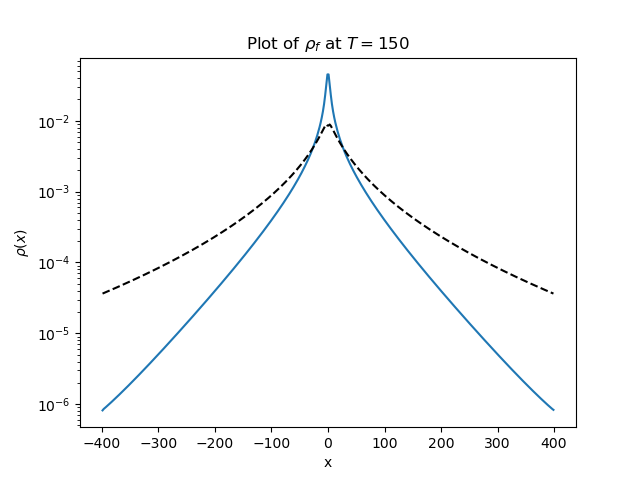}
  \end{subfigure}
  \begin{subfigure}{0.48\textwidth}
    \includegraphics[width=\linewidth]{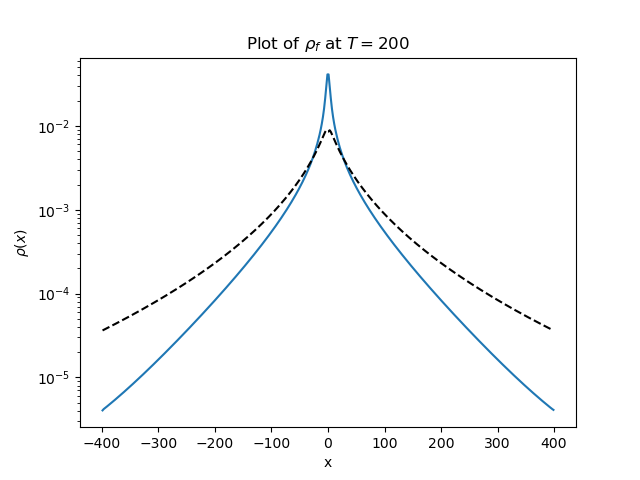}
  \end{subfigure}

  % Line 3
  \begin{subfigure}{0.48\textwidth}
    \includegraphics[width=\linewidth]{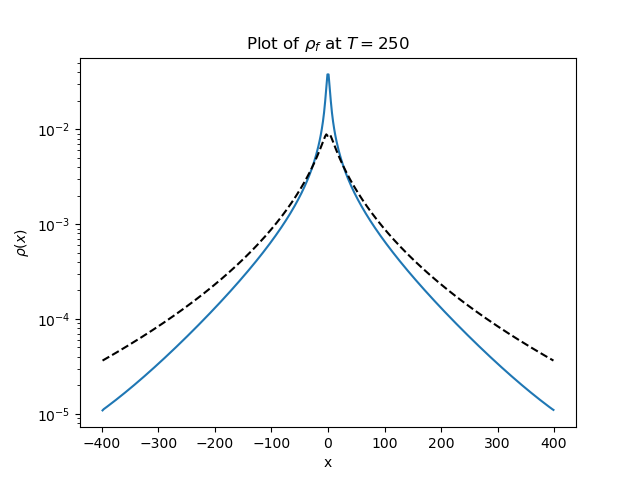}
  \end{subfigure}
  \begin{subfigure}{0.48\textwidth}
    \includegraphics[width=\linewidth]{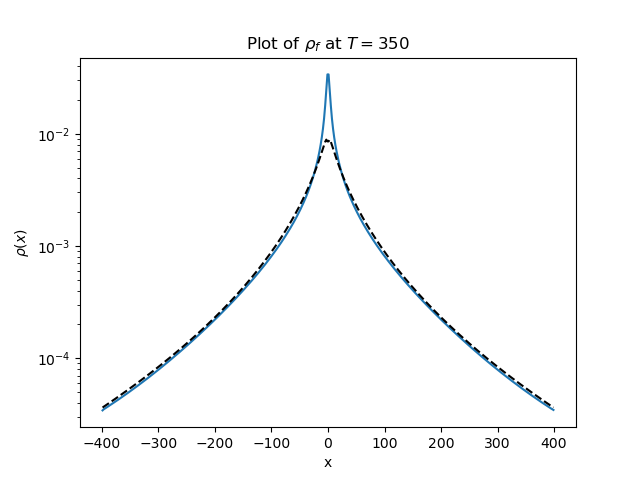}
  \end{subfigure}
  
  \caption{Blue line: evolution of the density $\rho_f$ of a solution $f$ to kinetic Fokker-Planck equation with $\alpha = 1.5$, $\beta=0.5$ and $\delta =1.15$. Black dashed line: expected asymptotic behaviour \eqref{eq:asympRhoG} with $\delta =1.15$.}
  \label{fig:alpha1.5_beta0.5}
\end{figure}

In a first example, we consider the choice of parameters $\alpha = 1.5$ and $\beta=0.5$, which falls under the assumptions of \Cref{thm:mainEXP}. We considered $L=400$ and $V=400$, the size of the grid of discretisation is $N_x=400$ and $N_v=400$. The final time is $T=350$ and the time step is $\Delta t=6.25\times 10^{-4}$, which is sufficiently small so that the CFL conditions are met on $[-L,L]\times[-V,V]$. The initial datum is the normalised function
\[
f_0(x,v) =\frac{1}{16} \exp\left(-\frac{\abs{x}}{2}-\frac{\abs{v}}{2}\right).
\]
In \Cref{fig:alpha1.5_beta0.5} we show the evolution of $\rho_f$. Notice that if we choose $\delta = 1.15$ in the definition of $\widetilde{G}$, we obtain that $\rho_{\widetilde{G}}$ is a good approximation of the tails of $\rho_f$. 

\begin{figure}
  \centering
  % Line 1
  \begin{subfigure}{0.48\textwidth}
    \includegraphics[width=\linewidth]{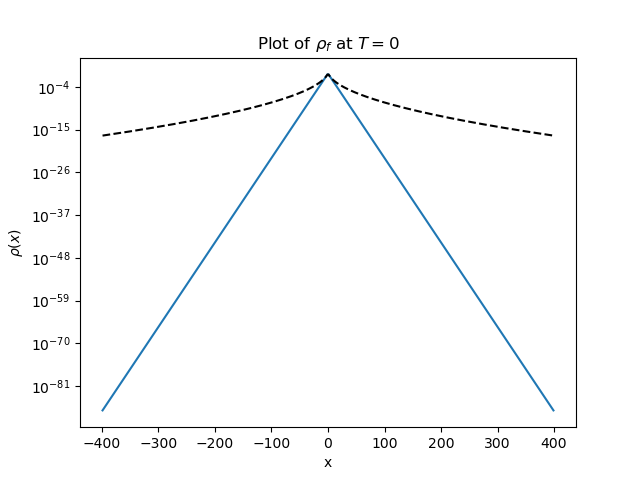}
  \end{subfigure}
  \begin{subfigure}{0.48\textwidth}
    \includegraphics[width=\linewidth]{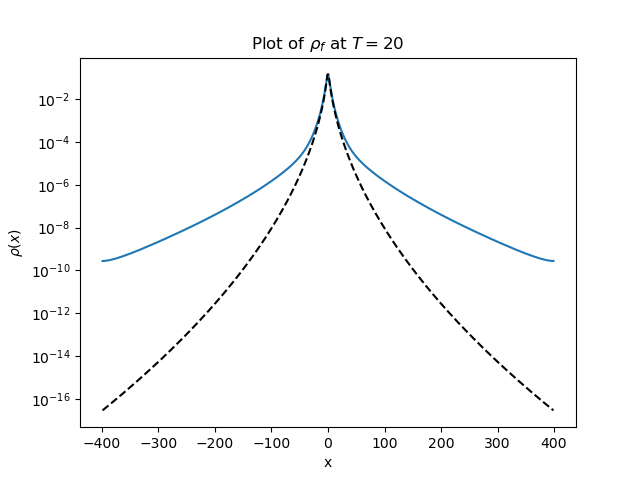}
  \end{subfigure}
  
  % Line 2
  \begin{subfigure}{0.48\textwidth}
    \includegraphics[width=\linewidth]{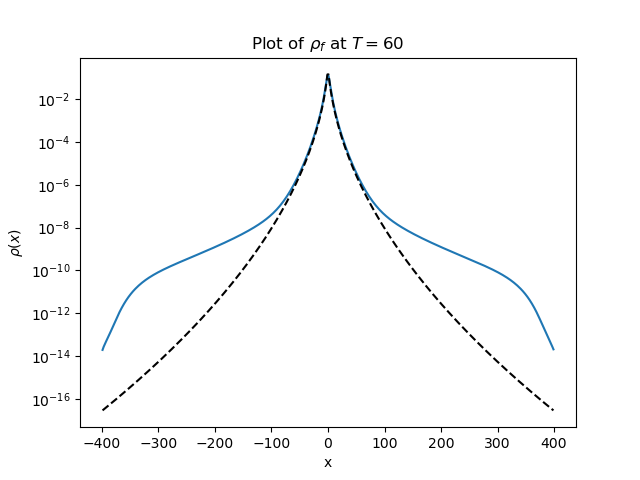}
  \end{subfigure}
  \begin{subfigure}{0.48\textwidth}
    \includegraphics[width=\linewidth]{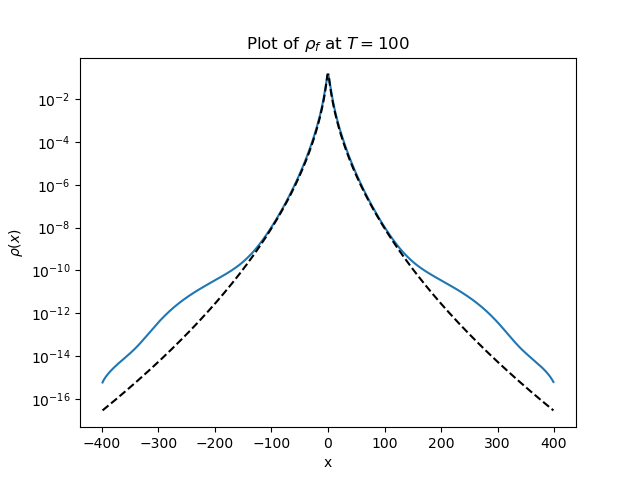}
  \end{subfigure}

  % Line 3
  \begin{subfigure}{0.48\textwidth}
    \includegraphics[width=\linewidth]{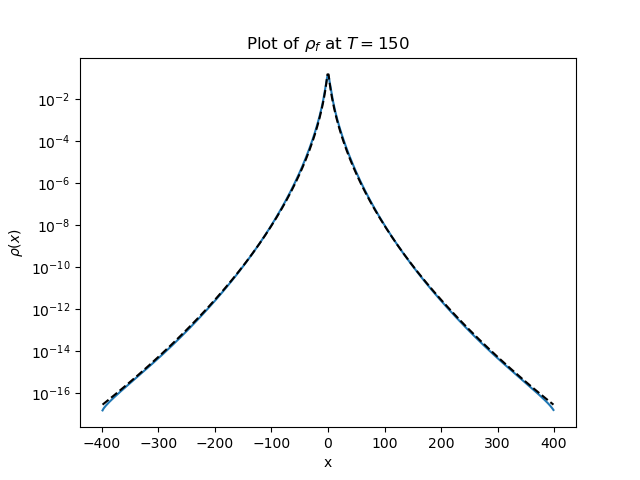}
  \end{subfigure}
  \begin{subfigure}{0.48\textwidth}
    \includegraphics[width=\linewidth]{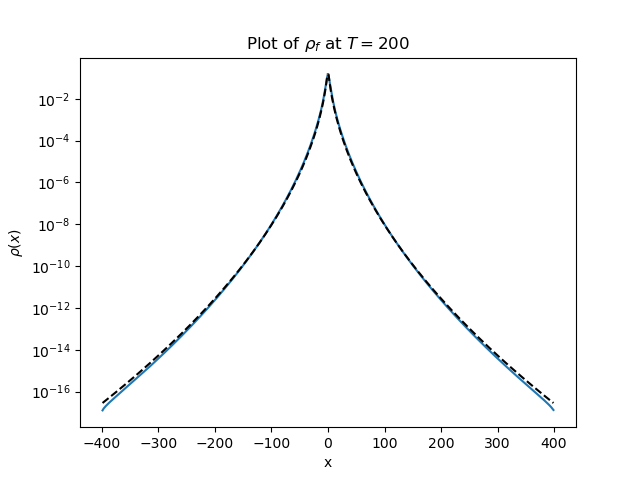}
  \end{subfigure}
  
  \caption{Blue line: evolution the density of the kinetic Fokker-Planck equation with $\alpha = 1$, $\beta=1$ and $\delta = 2$. Black dashed line: expected asymptotic behaviour \eqref{eq:asympRhoG} with $\delta=2$.}
  \label{fig:alpha1_beta_1}
\end{figure}

As an exploratory case, we also considered the case $\alpha = 1$ and $\beta = 1$, that is when the potential $V(x)$ grows linearly in space. This case is not covered by our theorem, but it is a limit case of the range of our admissible parameters. We considered the same domain and discretization of the previous example, that is $L,V= 400$, $N_x,N_v=400$, and the same initial datum. The final time is $T=200$ and in \Cref{fig:alpha1_beta_1}, we show the evolution of the density $\rho_f$. We observe a slightly different evolution of the density $\rho_f$, as the tails seems to spread more strongly, but after some time the solution stabilizes around a certain profile with the expected asymptotic behaviour \eqref{eq:asympRhoG}.

\subsection{Level plots of the steady state}
In this section we aim to compare the steady state $G$ of the kinetic Fokker-Planck equation and the function $\widetilde{G}(x,v) = \exp(-\delta E^{\frac{\beta}{2}})$ suggested by \Cref{thm:mainEXP}. We considered the final state $f(T,\cdot,\cdot)$ of the numerical simulation as an approximation of $G$.

\begin{figure}
  \centering
  % Line 1
  \begin{subfigure}{0.48\textwidth}
    \includegraphics[width=\linewidth]{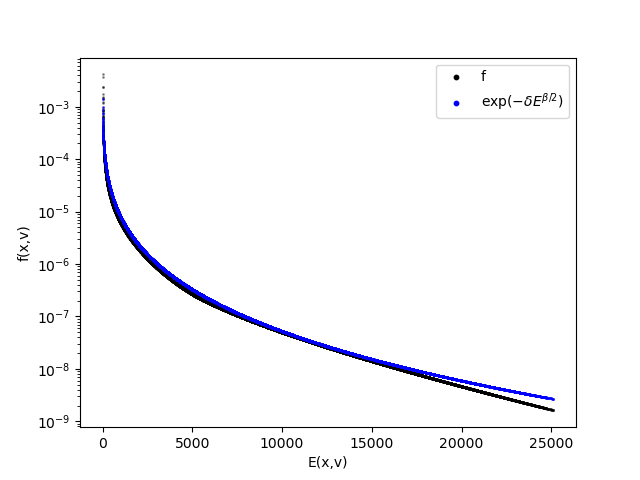}
  \end{subfigure}
  \begin{subfigure}{0.48\textwidth}
    \includegraphics[width=\linewidth]{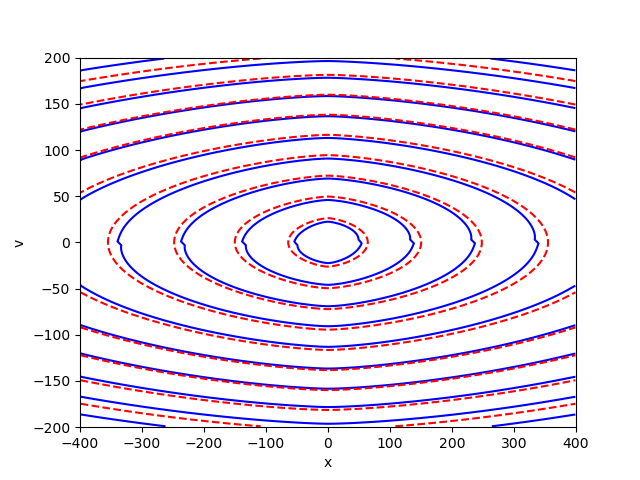}
  \end{subfigure}
  \caption{Profile and contour-plot of the solution $f$ of the kinetic Fokker-Planck equation at time $T=350$ for $\alpha = 1.5$ and $\beta=0.5$ and $\delta = 1.15$.}
  \label{fig:Energy_alpha1.5_beta0.5}
\end{figure}

At first, we would like to check whether the approximated steady state $f(T,x,v)$ resembles a function of the type $(x,v)\mapsto \Gamma(E(x,v))$ for a certain profile $\Gamma\colon\R_+\to\R_+$. To verify this, it is sufficient to check if all the points of the type $(E(x,v) , f(T,x,v))$ with $(x,v)\in\R\times\R$ fall approximately on the same curve of $\R^2$ (i.e. on a parametric curve in $\R^2$ of the form $s\mapsto(s,\Gamma(s)))$. In \Cref{fig:Energy_alpha1.5_beta0.5} we consider $\alpha=1.5$ and $\beta=0.5$. In the left picture, we show a scatter plot of all points of the type $(E(x,v),f(T,x,v))$ for all $(x,v)$ of the discretized domain $[-L,L]\times[-V,V]$. We observe that these points (in black) tend to gather around the same curve. Furthermore, by plotting the points of the type $(E(x,v), \exp(-\delta E(x,v)^{\frac{\beta}{2}}))$ in blue, we notice a close similarity with the cloud of points $(E(x,v), f(T,x,v))$. This is a sign that $G(x,v)$ is approximately a function of the energy, and this profile is $\Gamma(s)=\exp(-\delta s^{\frac{\beta}{2}})$. Further confirmation of this fact is given by the contour plot on the right of \Cref{fig:Energy_alpha1.5_beta0.5}. In that picture, we compare some level curves of $f(T,x,v)$ and $\exp(-\delta E^{\frac\beta2})$, we notice a very close similarity.

\begin{figure}
  \centering
  % Line 1
  \begin{subfigure}{0.48\textwidth}
    \includegraphics[width=\linewidth]{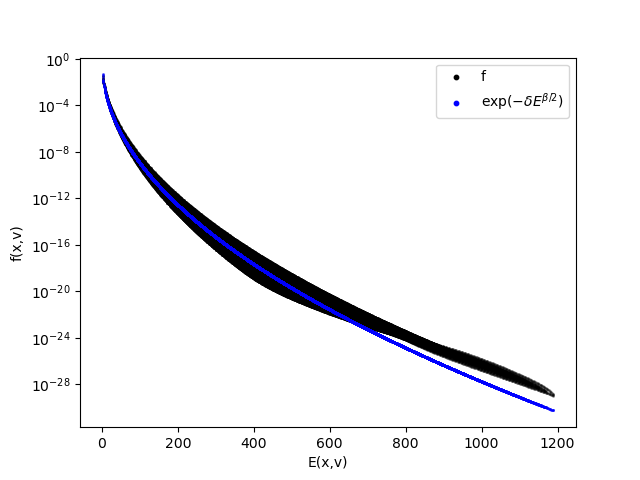}
  \end{subfigure}
  \begin{subfigure}{0.48\textwidth}
    \includegraphics[width=\linewidth]{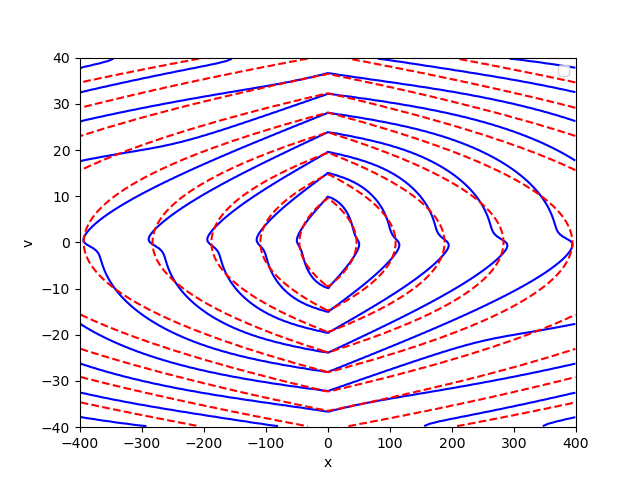}
  \end{subfigure}
  \caption{Profile and contour-plot of the steady state of the kinetic Fokker-Planck equation for $\alpha = 1$ and $\beta=1$ and $\delta = 2$.}
  \label{fig:Energy_alpha1_beta1}
\end{figure}

We carried out the same analysis for the case $\alpha=1$ and $\beta=1$, the results are shown in \Cref{fig:Energy_alpha1_beta1}. We can notice a greater discrepancy between the approximation of steady state $f(T,\cdot,\cdot)$ and the function $\exp(-\delta E^{\frac\beta2})$, both in the scatter plot on the left and in the contour plot on the right. In the left picture of \Cref{fig:Energy_alpha1_beta1}, the cloud of black points, represented by $(E(x,v),f(T,x,v))$, shows a much greater dispersion around a one-dimensional curve. This fact indicates some variability of the values of $f(T,\cdot,\cdot)$ over points $(x,v)$ with same energy $E$. This contrast is also visible from the contour plot on the right of \Cref{fig:Energy_alpha1_beta1}, in which we note that the level curves of $f(T,\cdot,\cdot)$ slightly deviate from the level curves of the energy. However, these oscillations appear quite restrained, so that they do not influence the decay of the density $\rho_f$ as we saw in the previous section. Indeed, apart from these small perturbations, the general trend of the steady state is very similar to $\exp(-\delta E^{\frac\beta2})$.

\section{Acknowledgments}
L. Z. is grateful to Richard Medina Rodrigues for fruitful discussions about the Harnack inequality and the De Giorgi method. L. Z. has received funding from the European Union’s Horizon 2020 research and innovation program under the Marie Skłodowska-Curie grant agreement No 945332.

\appendix
\section{Reminder on the Harris theorem }\label{sec:Harris}
Harris-type theorems are fundamental tools in the study of Markov processes and PDEs. Under suitable Lyapunov and positivity conditions, they ensure the existence of a stationary measure and convergence toward it. The seminal works for such results are due to Doeblin \cite{Doeblin1940} and Harris \cite{H56}. Doeblin established exponential convergence for Markov processes under the assumption of uniformly positive transition probabilities (the so-called Doeblin theorem \cite{CM21,EY23}), while Harris provided sufficient conditions for the existence of a unique stationary measure. Along the decades, these results have been refined to yield quantitative convergence estimates \cite{HM11,DMT95,MT09,MT93,CM21}; see also \cite{H16} for modern proofs of both theorems. Their development has been boosted as they have numerous applications to PDEs, where they complement hypocoercivity techniques for kinetic equations, see, for instance, Hu and Wang \cite{HW19}, Eberle et al. \cite{EGZ16}, Canizo et al. \cite{CCEY20}, Cao \cite{C21} and Laflèche \cite{L20}.

In this section we recall the statement of the sub-geometric Harris Theorem, which we took in particular from the more recent work \cite{CM21}. Let us recall the main notations. We denote by $\mathfrak{M}(\RRd)$ the set of all signed measures on $\RRd$ and  $\mathcal{P}(\RRd)\subseteq \mathfrak{M}(\RRd)$ the set of all probability measures. We consider the total variation norm
\[
\norm{\mu}_{TV} := \int_{\RRd} \d\abs{\mu}\qquad \text{ for every } \mu\in\mathfrak{M}(\RRd).
\]
For every positive weight function $m\colon\RRd\to [1,+\infty)$, we consider the weighted norm
\[
\norm{\mu}_{m} := \int_{\RRd} m\d\abs{\mu}\qquad \text{ for every } \mu\in\mathfrak{M}(\RRd).
\]
A stochastic semigroup is a family $(S(t))_{t\in[0,+\infty)}$ of mass preserving operators $S(t)\colon \mathcal{P}(\RRd)\to\mathcal{P}(\RRd)$ such that $S(0) = \mathsf{Id}_{\RRd}$ and $S(t)\circ S(s) =S(t+s) $ for all $s,t \geq 0$. When a stochastic semigroup is generated by an operator $\cL$, we denote the semigroup by $S_\cL$.

The main hypotheses of the Harris theorem are the following.
\begin{hypothesis}[Weak Lyapunov condition]\label{Hyp:Lyapunov}
There exists a continuous function $m\colon \R^d\times \R^d \to [1,+\infty)$ with pre-compact level sets such that 
\begin{equation}
\cL^* m\leq C\1_{B_R} -\phi(m) 
\end{equation}
for some constants $C,R>0$ and some strictly concave function $\phi\colon\R_+\to\R$ with $\phi(0)=0$ and increasing to infinity.
\end{hypothesis}

\begin{hypothesis}[Positivity condition]\label{Hyp:minorisation}
We say that the stochastic semigroup $S_\cL$ satisfies the positivity condition on a set $\mathcal{C}$ if there exists a probability measure $\mu_*$ and a constant $\eta\in(0,1)$  such that for a certain $T>0$
\begin{equation}\label{con:minorisation}
S_{\cL}(T)\mu \geq \eta\mu_*\int_\mathcal{C}\mu 
\end{equation}
for all positive measures $\mu$.
\end{hypothesis}

We now give the statement of the Harris theorem and we refer to \cite[Thm 5.6]{CM21} for the details of the proof.
\begin{thm}[Sub-geometric Harris Theorem] \label{thm:harris}
Consider a stochastic semigroup $S_\cL$ with generator $\cL$ that satisfies \Cref{Hyp:Lyapunov} for a continuous function $m\colon \RRd\to[1,\infty)$ and \Cref{Hyp:minorisation} in a set $\mathcal{C}=\{(x,v)\in\R^d\times\R^d\;\colon\; m(x,v)\leq M \}$ for large enough $M$. Then there exists a unique invariant measure $\mu_G\in \mathcal{P}(\R^d\times\R^d)$ such that
\begin{equation}\label{Harris:L1}
    \int_{\R^d\times\R^d}\phi(m)d\mu_G<\infty
\end{equation}
and there exists a decay-rate function $\Theta(t)$ such that
\begin{equation}\label{Harris:rate}
    \norm{S_{\cL}(t)\mu-\mu_G}_{TV}\lesssim \Theta(t)\norm{\mu-\mu_G}_m
\end{equation}
for any probability measure $\mu$. 
\end{thm}
The function $\Theta(t)$ can be explicitly computed from the concave function $\phi$ appearing in the Lyapunov condition, for all the details we refer to \cite[Sec. 4]{CM21}. 
In this work, $\phi$ can have, depending on the case, one of the two expressions we discuss here below.

The first case is $\phi(m)=m/(\log m)^\sigma$ for some $\sigma>0$. Following the examples of \cite[Sec. 4]{CM21}, we can compute that the function $\Theta(t)$ is given by
\[
\Theta(t)=e^{-\lambda t^{\frac{1}{1+\sigma}}}
\]
for an explicitly computable constant $\lambda>0$.

On the other hand, the function $\phi$ can take the expression $\phi(m)=m^{1-\kappa}$ for a certain $0<\kappa<1$. In this case, the decay function $\Theta(t)$ given from the Harris theorem is 
\[
\Theta(t)=\frac{1}{(1+t)^{\frac{1}{\kappa}}},
\]
which is a polynomial decay.

\bibliographystyle{plain} 
\bibliography{library} 

\begin{thebibliography}{10}

\bibitem{MR4017782}
F.~Anceschi, M.~Eleuteri, and S.~Polidoro.
\newblock A geometric statement of the {H}arnack inequality for a degenerate
  {K}olmogorov equation with rough coefficients.
\newblock {\em Communications in Contemporary Mathematics}, 21:1850057, 11
  2019.

\bibitem{Bakry2008}
D.~Bakry, P.~Cattiaux, and A.~Guillin.
\newblock Rate of convergence for ergodic continuous {M}arkov processes:
  Lyapunov versus poincaré.
\newblock {\em J. Funct. Anal.}, 254:727--759, 2008.

\bibitem{BCG26}
M.~Bisi, M.~Conte, and M.~Groppi.
\newblock Kinetic modeling of knowledge and wealth dynamics in national and
  global markets.
\newblock {\em Mathematical Models and Methods in Applied Sciences},
  36(06):1271--1305, 2026.

\bibitem{bony_1969}
J-M. Bony.
\newblock Principe du maximum, in\'egalit\'e de {Harnack} et unicit\'e du
  probl\`eme de {Cauchy} pour les op\'erateurs elliptiques d\'eg\'en\'er\'es.
\newblock {\em Annales de l'Institut Fourier}, 19(1):277--304, 1969.

\bibitem{BouinEmericandDolbeault2024}
E.~Bouin, J.~Dolbeault, and L.~Ziviani.
\newblock ${L}^2$-hypocoercivity methods for kinetic {F}okker-{P}lanck
  equations with factorised {G}ibbs states.
\newblock {\em Kolmogorov Operators and Their Applications}, pages 23--56,
  2024.

\bibitem{BEZ25}
E.~Bouin, J.~Evans, and L.~Ziviani.
\newblock Sub-exponential tails in biased run and tumble equations with
  unbounded velocities.
\newblock {\em arXiv:2505.08061}, 2025.

\bibitem{Cao2019}
C.~Cao.
\newblock The kinetic {F}okker-{P}lanck equation with weak confinement force.
\newblock {\em Commun. Math. Sci.}, 17:2281--2308, 2019.

\bibitem{C21}
C.~Cao.
\newblock The kinetic {F}okker–{P}lanck equation with general force.
\newblock {\em Journal of Evolution Equations}, 21:2293–2337, 06 2021.

\bibitem{CGMM24}
K.~Carrapatoso, P.~Gabriel, Medina R., and S.~Mischler.
\newblock Constructive {K}rein-{R}utman result for kinetic {F}okker-{P}lanck
  equations in a domain.
\newblock {\em arXiv:2407.10530}, 2024.

\bibitem{CCEY20}
J.~A. Cañizo, C.~Cao, J.~Evans, and H.~Yoldaş.
\newblock Hypocoercivity of linear kinetic equations via {H}arris's theorem.
\newblock {\em Kinet. Relat. Models}, 13(1):97--128, 2020.

\bibitem{CM21}
J.~A. Cañizo and S.~Mischler.
\newblock {H}arris-type results on geometric and subgeometric convergence to
  equilibrium for stochastic semigroups.
\newblock {\em J. Funct. Anal.}, 284:109830, 2023.

\bibitem{CC70}
J.S Chang and G~Cooper.
\newblock A practical difference scheme for fokker-planck equations.
\newblock {\em Journal of Computational Physics}, 6(1):1--16, 1970.

\bibitem{Cordier2005}
S.~Cordier, L.~Pareschi, and G.~Toscani.
\newblock On a kinetic model for a simple market economy.
\newblock {\em Journal of Statistical Physics}, 120(1):253--277, jul 2005.

\bibitem{DG56}
E.~De~Giorgi.
\newblock Sull'analiticit\`a{} delle estremali degli integrali multipli.
\newblock {\em Atti Accad. Naz. Lincei Rend. Cl. Sci. Fis. Mat. Nat. (8)},
  20:438--441, 1956.

\bibitem{DG57}
E.~De~Giorgi.
\newblock Sulla differenziabilit\`a{} e l'analiticit\`a{} delle estremali degli
  integrali multipli regolari.
\newblock {\em Mem. Accad. Sci. Torino. Cl. Sci. Fis. Mat. Nat. (3)}, 3:25--43,
  1957.

\bibitem{Doeblin1940}
W.~Doeblin.
\newblock Éléments d'une théorie générale des chaînes simples constantes
  de {M}arkoff.
\newblock {\em Annales scientifiques de l'École Normale Supérieure},
  57:61--111, 1940.

\bibitem{MR2576899}
J.~Dolbeault, C.~Mouhot, and C.~Schmeiser.
\newblock Hypocoercivity for kinetic equations with linear relaxation terms.
\newblock {\em C. R. Math. Acad. Sci. Paris}, 347(9-10):511--516, 2009.

\bibitem{Dolbeault2015}
J.~Dolbeault, C.~Mouhot, and C.~Schmeiser.
\newblock Hypocoercivity for linear kinetic equations conserving mass.
\newblock {\em Transactions of the American Mathematical Society},
  367:3807--3828, 2 2015.

\bibitem{MR2499863}
R.~Douc, G.~Fort, and A.~Guillin.
\newblock Subgeometric rates of convergence of {$f$}-ergodic strong {M}arkov
  processes.
\newblock {\em Stochastic Process. Appl.}, 119(3):897--923, 2009.

\bibitem{DMT95}
D.~Down, S.~P. Meyn, and R.~L. Tweedie.
\newblock Exponential and uniform ergodicity of {M}arkov processes.
\newblock {\em Ann. Probab.}, 23(4):1671--1691, 1995.

\bibitem{MR2813582}
R.~Duan.
\newblock Hypocoercivity of linear degenerately dissipative kinetic equations.
\newblock {\em Nonlinearity}, 24(8):2165--2189, 2011.

\bibitem{EGZ16}
A.~Eberle, A.~Guillin, and R.~Zimmer.
\newblock Quantitative harris type theorems for diffusions and mckean-vlasov
  processes.
\newblock {\em Transactions of the American Mathematical Society}, 371, 06
  2016.

\bibitem{EM24}
J.~Evans and A.~Menegaki.
\newblock Properties of non-equilibrium steady states for the non-linear {BGK}
  equation on the torus.
\newblock {\em Ann. Inst. H. Poincaré C Anal. Non Linéaire (2024), published
  online first}, 2024.

\bibitem{EY23}
J.~Evans and H.~Yoldaş.
\newblock On the asymptotic behavior of a run and tumble equation for bacterial
  chemotaxis.
\newblock {\em SIAM Journal of Mathematical Analysis}, 55(6), 2023.

\bibitem{EY24}
J.~Evans and H~Yoldaş.
\newblock Trend to equilibrium for run and tumble equations with non-uniform
  tumbling kernels.
\newblock {\em Acta Applicandae Mathematicae}, 191(6), 2024.

\bibitem{MR3923847}
F.~Golse, C.~Imbert, C.~Mouhot, and A.~F. Vasseur.
\newblock Harnack inequality for kinetic {F}okker-{P}lanck equations with rough
  coefficients and application to the {L}andau equation.
\newblock {\em Ann. Sc. Norm. Super. Pisa Cl. Sci. (5)}, 19:253--295, 2019.

\bibitem{Guerand_Imbert_2023}
J.~Guerand and C.~Imbert.
\newblock Log-transform and the weak {H}arnack inequality for kinetic
  {F}okker-{P}lanck equations.
\newblock {\em Journal of the Institute of Mathematics of Jussieu},
  22(6):2749–2774, 2023.

\bibitem{GM22}
J.~Guerand and C.~Mouhot.
\newblock Quantitative {D}e {G}iorgi methods in kinetic theory.
\newblock {\em Journal de l’École polytechnique - Mathématiques},
  9:1159--1181, 2022.

\bibitem{H16}
M.~Hairer.
\newblock Convergence of {M}arkov processes.
\newblock {\em Unpublished lecture notes}, January 2016.

\bibitem{HM11}
M.~Hairer and J.~C. Mattingly.
\newblock Yet another look at {H}arris' ergodic theorem for {M}arkov chains.
\newblock In Robert Dalang, Marco Dozzi, and Francesco Russo, editors, {\em
  Seminar on Stochastic Analysis, Random Fields and Applications VI}, pages
  109--117, Basel, 2011. Springer Basel.

\bibitem{H56}
T.~E. {H}arris.
\newblock The existence of stationary measures for certain {M}arkov processes.
\newblock In {\em Proceedings of the {T}hird {B}erkeley {S}ymposium on
  {M}athematical {S}tatistics and {P}robability, 1954--1955, vol. {II}}, pages
  113--124. University of California Press, Berkeley and Los Angeles, 1956.

\bibitem{Helffer2005}
B.~Helffer and F.~Nier.
\newblock {\em Hypoelliptic estimates and spectral theory for {F}okker-{P}lanck
  operators and Witten Laplacians}, volume 1862.
\newblock Springer-Verlag, Berlin, 2005.

\bibitem{MR2294477}
F.~H\'erau.
\newblock Short and long time behavior of the {F}okker-{P}lanck equation in a
  confining potential and applications.
\newblock {\em J. Funct. Anal.}, 244(1):95--118, 2007.

\bibitem{HW19}
S.~Hu and X.~Wang.
\newblock Subexponential decay in kinetic {F}okker–{P}lanck equation: Weak
  hypocoercivity.
\newblock {\em Bernoulli}, 25:174--188, 02 2019.

\bibitem{Hrau2004}
F.~Hérau and F.~Nier.
\newblock Isotropic hypoellipticity and trend to equilibrium for the
  {F}okker-{P}lanck equation with a high-degree potential.
\newblock {\em Arch. Ration. Mech. Anal.}, 171:151--218, 2004.

\bibitem{MR4049224}
C.~Imbert and L.~Silvestre.
\newblock The weak {H}arnack inequality for the {B}oltzmann equation without
  cut-off.
\newblock {\em J. Eur. Math. Soc. (JEMS)}, 22(2):507--592, 2020.

\bibitem{Kavian2021}
O.~Kavian, S.~Mischler, and M.~Ndao.
\newblock The {F}okker-{P}lanck equation with subcritical confinement force.
\newblock {\em J. Math. Pures Appl. (9)}, 151:171--211, 2021.

\bibitem{KT00}
A.~Kurganov and E.~Tadmor.
\newblock New high-resolution central schemes for nonlinear conservation laws
  and convection–diffusion equations.
\newblock {\em Journal of Computational Physics}, 160(1):241--282, 2000.

\bibitem{L20}
L.~Lafleche.
\newblock Fractional {F}okker--{P}lanck equation with general confinement
  force.
\newblock {\em SIAM Journal on Mathematical Analysis}, 52(1):164--196, 2020.

\bibitem{Lanconelli_1994}
E.~Lanconelli and S.~Polidoro.
\newblock On a class of hypoelliptic evolution operators.
\newblock {\em Rend. Semin. Mat., Torino}, 52(1):29--63, 1994.

\bibitem{LRR22}
T.~Lelièvre, M.~Ramil, and J.~Reygner.
\newblock A probabilistic study of the kinetic {F}okker–{P}lanck equation in
  cylindrical domains.
\newblock {\em Journal of Evolution Equations}, 22:38, 2022.

\bibitem{MT93}
S.~P. Meyn and R.~L. Tweedie.
\newblock {Stability of {M}arkovian processes III: Foster-Lyapunov criteria for
  continuous-time processes.}
\newblock {\em Adv. Appl. Prob.}, pages 518--548, 1993.

\bibitem{MT09}
S.~P. Meyn and R.~L. Tweedie.
\newblock {\em {M}arkov {C}hains and {S}tochastic {S}tability}.
\newblock 2nd ed., Cambridge Mathematical Library. Cambridge: Cambridge
  University Press, 2019.
\newblock With a prologue by P. W. Glynn.

\bibitem{PP04}
A.~Pascucci and S.~Polidoro.
\newblock The {M}oser's iterative method for a class of ultraparabolic
  equations.
\newblock {\em Communications in Contemporary Mathematics}, 06:395--417, 6
  2004.

\bibitem{RB24}
S.~Roy and A.~Borzì.
\newblock Numerical approximation of kinetic fokker–planck equations with
  specular reflection boundary conditions.
\newblock {\em Journal of Computational Physics}, 503:112841, 2024.

\bibitem{T14}
N.M. Temme.
\newblock {\em Asymptotic Methods for Integrals}.
\newblock World Scientific, 12 2014.

\bibitem{T06}
G.~Toscani.
\newblock {Kinetic models of opinion formation}.
\newblock {\em Communications in Mathematical Sciences}, 4(3):481 -- 496, 2006.

\bibitem{Villani2009}
C.~Villani.
\newblock Hypocoercivity.
\newblock {\em Mem. Amer. Math. Soc.}, 202:iv+141, 2009.

\bibitem{Yolda2023}
H~Yoldaş.
\newblock On quantitative hypocoercivity estimates based on {H}arris-type
  theorems.
\newblock {\em Journal of Mathematical Physics}, 64:031101--1 -- 031101--28,
  2023.

\end{thebibliography}

\end{document}